\documentclass[11pt,a4paper]{article}
\usepackage[reqno,sumlimits]{amsmath}
\usepackage{latexsym}
\usepackage{graphics}
\usepackage{amsfonts} 
\usepackage{amssymb}   
\usepackage{theorem}
\usepackage{epic,eepic}
\usepackage{epsfig}
\usepackage{amscd}
\newtheorem{theorem}{Theorem}[section]
\newtheorem{proposition}[theorem]{Proposition}

\newtheorem{lemma}[theorem]{Lemma}

{\theoremstyle{plain}
{\theorembodyfont{\normalfont\rmfamily}
\newtheorem{definition}[theorem]{Definition} 
\newtheorem{remark}[theorem]{Remark}

}}
\newcommand{\cqd}{\hfill $\Box$\medskip}
\newcommand{\proof}{\noindent{\bf Proof: }}
\newcommand{\abs}[1]{\left| #1\right|}
\newcommand{\norm}[1]{\left\| #1\right\|}

\newcommand{\id} {\mbox{\rm Id}}
\newcommand{\im} {\mbox{\rm Im\,}}

\newcommand{\Span}{\mbox{\rm span}}

\newcommand{\A}{{\cal A}}
\newcommand{\BB}{{\cal B}}
\newcommand{\CC}{{\cal C}}
\newcommand{\DD}{{\cal D}}
\newcommand{\EE}{{\cal E}}
\newcommand{\FF}{{\cal F}}

\newcommand{\II}{{\cal I}}
\newcommand{\JJ}{{\cal J}}

\newcommand{\LL}{{\cal L}}

\newcommand{\OO}{{\cal O}}

\newcommand{\RR}{{\cal R}}
\newcommand{\TT}{{\cal T}}
\newcommand{\UU}{{\cal U}}

\newcommand{\WW}{{\cal W}}

\newcommand{\Cc}{{\mathbb{C}}}

\newcommand{\Ee}{{\mathbb{E}}}

\newcommand{\Ii}{{\mathbb{I}}}

\newcommand{\Pp}{{\mathbb{P}}}
\newcommand{\Qq}{{\mathbb{Q}}}
\newcommand{\Rr}{{\mathbb{R}}}
\newcommand{\Ss}{{\mathbb{S}}}
\newcommand{\Tt}{{\mathbb{T}}}

\newcommand{\Zz}{{\mathbb{Z}}}
\begin{document}
\begin{titlepage}

\title
{\bf Renormalisation of flows on the multidimensional torus close to a
$KT$ frequency vector}

\author{Jo\~ao Lopes Dias\thanks{E-mail: j.lopes-dias@damtp.cam.ac.uk}\\   
{\small Department of Applied Mathematics and Theoretical Physics,} \\
{\small University of Cambridge} \\ 
{\small Silver Street, Cambridge  CB3 9EW, England}}
\date{April, 2001}

\maketitle
\begin{abstract}
We use a renormalisation operator $\RR$ acting on a space of vector
fields on $\Tt^d$, $d\geq2$, to prove the existence of a submanifold
of vector fields equivalent to constant.
The result comes from the existence of a fixed point $\omega$ of
$\RR$ which is hyperbolic.
This is done for a certain class $KT_d$ of frequency vectors
$\omega\in\Rr^d$, called of Koch type.
The transformation $\RR$ is constructed using a time rescaling, a
linear change of basis plus a periodic non-linear map
isotopic to the identity, which we derive by a ``homotopy trick''.
\end{abstract}

\end{titlepage}

%%%%%%%%%%%%%%%%%%%%%%%%%%%%%%%%%%%%%%%%%%%%%%%%%%%%%%%%%%%%%%%%%%%%%%%%%%%%

\section{Introduction}

In this paper we consider a renormalisation transformation $\RR$ that
acts on a Banach space of analytic vector fields $X$ on $\Tt^d=\Rr^d/\Zz^d$, $d\geq2$.
Its domain is an open ball  
around a non-zero constant vector field
$\omega\in\Rr^d$ (that generates a linear flow).
The operator $\RR$ is a $C^1$ infinite-dimensional dynamical system
with a hyperbolic fixed point at $\omega$.
We show that the stable and unstable manifolds, given by the local dynamics of $\RR$
around $\omega$, yield different equivalence classes of vector fields.
The elements of the codimension-$(d-1)$ stable manifold are equivalent
to $\omega$,
as $\RR$ is designed to be an equivalence between flows.
In other words, the renormalisation asymptotically contracts locally any $X$ into a
$(d-1)$-parameter family of constant vector fields that includes $\omega$.
The map $\RR$ basically consists of the composition of a time
rescaling with a diffeomorphism of $\Tt^d$ isotopic to a linear
automorphism in $GL(d,\Zz)$.

This idea of renormalising vector fields is due to MacKay
\cite{MacKay}, and
we follow an approach inspired by the
work of Koch \cite{Koch} on
renormalisation for $d$-degrees of freedom analytic Hamiltonian systems.
The latter applies to the problem of stability of invariant tori associated with a frequency
vector $\omega$.

The linear transformation is the main feature of the operator.
It is a change of the basis of
$\Tt^d$ lifted to $\Rr^d$, made to
enlarge the region around the orbits of the linear flow.
By iterating this, we look closer at chosen regions, but
periodicity and the whole torus are kept at each stage.

We also need to rescale the time as orbits take longer to cross the
``new'' torus.
The topological characteristics of the trajectories are not affected
by such transformation.

The non-linear part of $\RR$, isotopic to the identity, is obtained by a
homotopy method to reduce the perturbation.
We use a flow of coordinate changes (diffeomorphisms) to find one which fully
eliminates some of the Fourier terms of $X$.
Those are chosen to be the easiest to do so, since we keep the
``resonant'' terms associated with ``small denominators'' (as it
appears in the usual KAM theory).
We are not concerned with trying to cancel all the perturbation terms
because the linear transformation shifts resonant
to ``non-resonant'' terms.
Eventually all perturbation terms are eliminated while iterating $\RR$.
The separation between resonant and non-resonant terms appears in \cite{Koch}.

If a vector field is a fixed point of the above procedure, it means that
its orbits exhibit self-similarity between the different scales. 
This is trivially deduced for the fixed point $\omega$.

The existence of the linear part, based on an arithmetical
condition \cite{Koch}, depends on certain conditions
imposed on $\omega$.
In particular, for $d=2$, that corresponds to the
set of vectors with a ``quadratic irrational'' slope.
Specifically for the two-dimensional case, it is defined in \cite{jld2}
a family of renormalisation iterative schemes
allowing a full Lebesgue measure set of diophantine vectors.
The methods involved therein are not easily generalisable to
higher dimensions, it is required a suitable choice of a multidimensional continued
fraction expansion algorithm.

\medskip

This paper is organised as follows.
In Section \ref{Equivalence of flows} we recall the basic ingredients
of equivalence of flows, and in
Section \ref{Motivating example for d=2} we introduce the
renormalisation idea through an example.
We rigorously construct the operator $\RR$ in Sections \ref
{Space of Analytic Vector Fields}, \ref{section
elimination far from resonance terms for local Td},
\ref{Change of Basis} and \ref {section renorm oper}, with proofs
of the statements also given in Section \ref{Proofs of Theorems}.
Finally, in Section \ref{Hyperbolicity of the Fixed Point omega} 
we present the main result using the
spectral properties of the derivative at the fixed point, with 
further discussion in Section \ref{Spectral
Properties and Invariant Manifolds}.

%%%%%%%%%%%%%%%%%%%%%%%%%%%%%%%%%%%%%%%%%%%%%%%%%%%%%%%%%%%%%%%%%%%%%%%%%%%%

\section{Equivalence of flows}\label{Equivalence of flows}

Consider a continuous vector field
$X$ on $\Tt^d$ that generates the flow $\phi_t$,
$\dot\theta=X(\theta)$, $\theta\in \Tt^d$, with lift
$\Phi_t$ to the universal cover $\Rr^d$.
Choosing a norm $\|\cdot\|$ in $\Rr^d$, we define 
$w_X(\theta_0)=\lim_{t\to\infty} \Phi_t(\theta_0)/
\|\Phi_t(\theta_0)\|$
to be the {\it winding ratio} of $X$ for the orbit of
$\theta_0\in\Tt^d$, if the limit exists and $\lim_{t\to\infty}\|\Phi_t(\theta_0)\|=\infty$.
Otherwise, if $\Phi_t(\theta_0)$ is bounded we put $w_X(\theta_0)=0$.
Also, if the limit does not exist or if
$\|\Phi_t(\theta_0)\|$ is unbounded but does not
tend to infinity, the winding ratio is not defined.

We say that two flows $\phi_t,\psi_t\colon \Tt^d\to\Tt^d$ are {\it
$C^r$-equivalent} if there is a $C^r$-diffeomorphism
$h\colon\Tt^d\to\Tt^d$ taking orbits of $\phi_t$ onto those of
$\psi_t$, preserving orientation. 
We are allowing to have
$h(\phi_t(\theta))=\psi_{\tau(h(\theta),t)} (h(\theta))$, where
$\tau(\theta,\cdot)$ is a homeomorphism of $\Rr$ for any $\theta\in\Tt^d$.
This is the same to say that two
vector fields $X$, $Y$ on $\Tt^d$ are equivalent if $\tau'X\circ h=Dh\,Y$
where $h$ is a flow equivalence as above and $\tau'$ the time derivative of $\tau$.
We emphasise that this relaxation of the usual requirement that $t$ be
preserved ($h\circ\phi_t=\psi_t\circ h$) provides more satisfactory equivalence
classes for flows, since we are in fact mainly interested in qualitative (topological)
properties of the flow.
Every $C^r$-equivalence is isotopic to a map with a lift to $\Rr^d$ in the group of the
linear automorphisms of the lattice $\Zz^d$ with determinant $\pm1$,
$GL(d,\Zz)$. 
The isotopy is given by periodic homeomorphisms on $\Tt^d$.
These kind of coordinate changes preserve the $d$-torus structure and
volume.

The set of
winding ratios of a flow on $\Tt^d$ generated by $X$ is
$w_X= \{w_X(\theta)\colon \theta\in\Tt^d\}$,
which is called the {\it
winding set}.
Any automorphism $T\colon \Rr^d\to\Rr^d$ induces the map $\hat
T\colon \Ss^{d-1}\to \Ss^{d-1}$ given by $x\mapsto Tx/\|Tx\|$, where
$\Ss^{d-1}=\{x\in\Rr^d\colon \|x\|=1\}$.
Let $T$ be in $GL(d,\Zz)$.
Then the winding set 
is invariant up to the action of $\hat T$, with a $C^r$-equivalence
$h$ isotopic to a map with lift $T$.
That is, if $X'=(Dh)^{-1}\,X\circ h$, then $w_{X'}(\theta')= \hat T
w_X(\theta)$, where $\theta'=h^{-1}(\theta)$.
In particular, the winding set is preserved by transformations isotopic to
the identity.

%%%%%%%%%%%%%%%%%%%%%%%%%%%%%%%%%%%%%%%%%%%%%%%%%%%%%%%%%%%%%%%%%%%%%%%%%%%%

\section{Motivating example for $d=2$}\label{Motivating example for d=2}

We start by motivating the renormalisation procedure with
a simple example.
Consider the case $d=2$ and the linear flow
described by the differential equation
$\dot \theta=\omega$,
with $\omega\in\Rr^2$ and $\theta\in\Tt^2$.
Given an initial condition $\theta_0\in\Tt^2$ for $t=0$, the solution of the
flow is 
$\phi_t(\theta_0)=\theta_0+\omega t \mod 1$, $t>0$.
The motion is simply a rotation with frequency vector $\omega$ and
winding ratio $w_{\omega}=\omega/\|\omega\|$. 
If, for all non-zero integer vectors $k$, $k\cdot\omega:=k_1\omega_1 +k_2\omega_2 \not=0$
(i.e. the slope $\omega_2/\omega_1$ is irrational, assuming $\omega_1$
does not vanish), then all the orbits are dense in
the torus (the flow is minimal).
Otherwise, they are closed curves (periodic orbits).

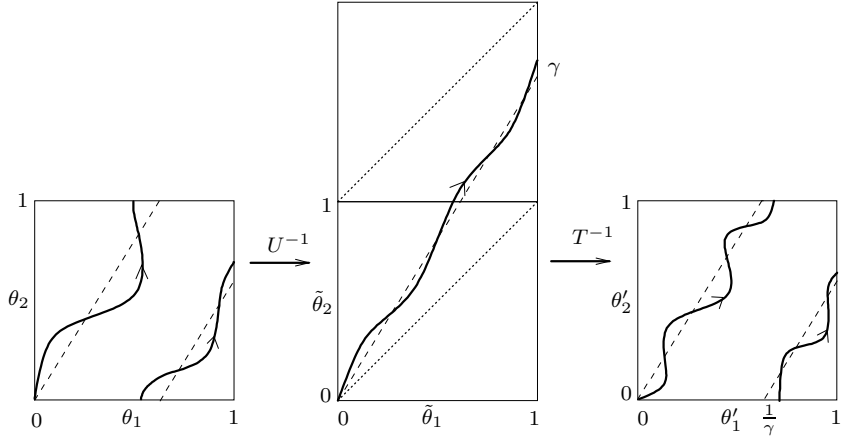
\begin{figure}
\begin{center}
\setlength{\unitlength}{0.00058333in}
\begingroup\makeatletter\ifx\SetFigFont\undefined%
\gdef\SetFigFont#1#2#3#4#5{%
  \reset@font\fontsize{#1}{#2pt}%
  \fontfamily{#3}\fontseries{#4}\fontshape{#5}%
  \selectfont}%
\fi\endgroup%
{\renewcommand{\dashlinestretch}{30}
\begin{picture}(8207,3894)(0,-10)
\thicklines
\path(4904,1528)(5429,1528)
\thinlines
\path(5309.000,1498.000)(5429.000,1528.000)(5309.000,1558.000)
\put(5091,1678){\makebox(0,0)[lb]{\smash{{{\SetFigFont{8}{9.6}{\familydefault}{\mddefault}{\updefault}$T^{-1}$}}}}}
\thicklines
\path(2189,1513)(2714,1513)
\thinlines
\path(2594.000,1483.000)(2714.000,1513.000)(2594.000,1543.000)
\put(2339,1611){\makebox(0,0)[lb]{\smash{{{\SetFigFont{8}{9.6}{\familydefault}{\mddefault}{\updefault}$U^{-1}$}}}}}
\put(4874,3237){\makebox(0,0)[lb]{\smash{{{\SetFigFont{8}{9.6}{\rmdefault}{\mddefault}{\updefault}$\gamma$}}}}}
\dashline{60.000}(240,269)(1365,2069)
\dashline{60.000}(5682,274)(6807,2074)
\drawline(6807,2074)(6807,2074)
\path(5682,2074)(7482,2074)(7482,274)
	(5682,274)(5682,2074)
\path(2977,3867)(4777,3867)(4777,2067)
	(2977,2067)(2977,3867)
\dashline{60.000}(2977,267)(4777,3192)
\dottedline{45}(2977,267)(4777,2067)
\dottedline{45}(2977,2067)(4777,3867)
\path(3989,2195)(4117,2247)(4124,2120)
\path(2977,2067)(4777,2067)(4777,267)
	(2977,267)(2977,2067)
\path(1154,1370)(1206,1490)(1259,1370)
\path(6351,1205)(6456,1198)(6419,1100)
\dashline{60.000}(1371,265)(2024,1325)
\dashline{60.000}(6831,288)(7484,1348)
\path(1799,785)(1860,845)(1897,740)
\path(239,2075)(2039,2075)(2039,275)
	(239,275)(239,2075)
\path(7304,860)(7394,905)(7402,793)
\thicklines
\path(2977,267)(2978,270)(2980,276)
	(2985,288)(2991,306)(3000,329)
	(3010,359)(3023,392)(3037,430)
	(3052,469)(3068,509)(3084,548)
	(3099,586)(3115,622)(3129,656)
	(3144,688)(3158,717)(3171,743)
	(3184,768)(3198,790)(3211,811)
	(3224,831)(3238,849)(3252,867)
	(3267,884)(3282,900)(3297,916)
	(3314,932)(3331,948)(3349,964)
	(3367,980)(3386,996)(3406,1013)
	(3426,1030)(3446,1047)(3466,1065)
	(3487,1083)(3507,1102)(3527,1120)
	(3547,1140)(3567,1160)(3586,1180)
	(3605,1201)(3623,1223)(3640,1245)
	(3657,1268)(3674,1292)(3690,1317)
	(3702,1338)(3714,1360)(3727,1383)
	(3739,1407)(3751,1432)(3763,1459)
	(3776,1486)(3788,1514)(3801,1543)
	(3813,1573)(3826,1604)(3839,1635)
	(3851,1667)(3864,1699)(3877,1731)
	(3890,1764)(3903,1796)(3915,1828)
	(3928,1859)(3941,1891)(3953,1921)
	(3966,1951)(3978,1979)(3991,2007)
	(4003,2034)(4015,2060)(4027,2085)
	(4040,2109)(4052,2132)(4065,2154)
	(4080,2181)(4097,2207)(4114,2232)
	(4131,2256)(4149,2279)(4168,2302)
	(4187,2325)(4207,2347)(4227,2369)
	(4247,2391)(4267,2412)(4288,2433)
	(4308,2454)(4328,2474)(4348,2494)
	(4368,2514)(4387,2533)(4405,2552)
	(4423,2571)(4440,2590)(4457,2609)
	(4472,2628)(4487,2647)(4502,2667)
	(4515,2686)(4528,2705)(4540,2726)
	(4552,2748)(4564,2771)(4577,2796)
	(4589,2822)(4602,2851)(4615,2882)
	(4628,2915)(4642,2951)(4656,2989)
	(4670,3028)(4685,3069)(4699,3110)
	(4713,3151)(4727,3190)(4739,3226)
	(4750,3258)(4759,3286)(4766,3307)
	(4771,3323)(4774,3334)(4776,3339)(4777,3342)
\path(239,283)(239,286)(240,293)
	(242,305)(245,324)(249,349)
	(255,380)(261,415)(269,455)
	(277,496)(287,537)(297,579)
	(307,618)(318,656)(330,691)
	(342,723)(354,752)(368,779)
	(382,804)(397,826)(414,846)
	(431,865)(450,882)(471,898)
	(492,912)(513,925)(537,938)
	(561,950)(587,962)(614,973)
	(643,985)(672,996)(702,1007)
	(733,1019)(764,1030)(796,1041)
	(828,1053)(859,1064)(890,1076)
	(920,1088)(949,1100)(977,1113)
	(1004,1126)(1029,1139)(1053,1153)
	(1075,1167)(1095,1182)(1113,1198)
	(1130,1214)(1145,1232)(1160,1252)
	(1172,1273)(1183,1296)(1192,1320)
	(1199,1345)(1205,1372)(1210,1400)
	(1213,1430)(1215,1460)(1215,1490)
	(1214,1522)(1213,1553)(1210,1585)
	(1206,1616)(1202,1646)(1197,1676)
	(1192,1705)(1187,1733)(1181,1759)
	(1176,1784)(1170,1808)(1165,1830)
	(1160,1851)(1155,1870)(1148,1900)
	(1142,1926)(1138,1949)(1134,1969)
	(1132,1988)(1130,2004)(1129,2020)
	(1129,2034)(1129,2046)(1130,2056)
	(1130,2062)(1131,2066)(1131,2068)
\path(5696,282)(5699,283)(5706,286)
	(5719,291)(5736,297)(5757,306)
	(5781,317)(5805,329)(5830,342)
	(5852,355)(5872,368)(5890,381)
	(5904,395)(5916,410)(5925,425)
	(5932,441)(5937,458)(5940,477)
	(5941,496)(5940,516)(5939,537)
	(5937,560)(5934,585)(5931,610)
	(5928,637)(5925,664)(5922,692)
	(5919,721)(5918,749)(5917,777)
	(5918,804)(5920,830)(5924,855)
	(5930,879)(5937,901)(5947,922)
	(5959,941)(5974,959)(5988,973)
	(6004,986)(6022,998)(6041,1010)
	(6063,1021)(6086,1032)(6110,1043)
	(6136,1053)(6163,1064)(6190,1074)
	(6218,1084)(6246,1094)(6274,1104)
	(6302,1114)(6329,1125)(6354,1135)
	(6379,1146)(6402,1157)(6423,1168)
	(6443,1180)(6460,1192)(6476,1205)
	(6489,1219)(6500,1233)(6510,1250)
	(6518,1268)(6523,1287)(6526,1308)
	(6528,1330)(6528,1353)(6526,1378)
	(6523,1403)(6519,1429)(6514,1455)
	(6508,1481)(6502,1508)(6495,1534)
	(6489,1559)(6483,1584)(6478,1608)
	(6474,1630)(6470,1651)(6468,1671)
	(6467,1689)(6468,1706)(6470,1722)
	(6476,1739)(6484,1754)(6494,1768)
	(6508,1781)(6524,1793)(6542,1803)
	(6561,1812)(6583,1820)(6606,1828)
	(6629,1835)(6653,1841)(6676,1846)
	(6699,1851)(6721,1856)(6742,1861)
	(6762,1866)(6780,1871)(6796,1876)
	(6813,1883)(6827,1892)(6840,1901)
	(6852,1913)(6862,1926)(6871,1941)
	(6880,1958)(6888,1978)(6895,1998)
	(6901,2018)(6906,2037)(6909,2053)
	(6912,2065)(6913,2072)(6914,2075)
\path(6964,266)(6964,267)(6964,270)
	(6964,278)(6964,293)(6965,315)
	(6965,343)(6967,375)(6968,411)
	(6970,448)(6973,484)(6976,519)
	(6979,550)(6984,579)(6989,605)
	(6994,628)(7001,648)(7009,665)
	(7017,680)(7027,693)(7038,704)
	(7052,715)(7067,725)(7084,734)
	(7102,741)(7122,748)(7142,754)
	(7164,760)(7186,765)(7208,771)
	(7230,776)(7252,782)(7273,788)
	(7292,795)(7310,803)(7327,812)
	(7342,821)(7355,832)(7366,845)
	(7375,857)(7382,872)(7389,888)
	(7394,905)(7398,924)(7401,945)
	(7403,966)(7405,989)(7406,1012)
	(7406,1036)(7405,1060)(7405,1085)
	(7404,1108)(7402,1131)(7401,1154)
	(7400,1175)(7399,1195)(7398,1214)
	(7398,1231)(7398,1248)(7399,1273)
	(7402,1295)(7407,1315)(7414,1333)
	(7424,1351)(7436,1369)(7449,1386)
	(7462,1401)(7473,1412)(7480,1419)
	(7483,1423)(7484,1423)
\path(1198,279)(1198,280)(1199,283)
	(1200,291)(1203,302)(1207,315)
	(1211,327)(1216,339)(1221,350)
	(1228,361)(1237,372)(1243,380)
	(1250,388)(1258,396)(1267,406)
	(1277,415)(1288,425)(1301,435)
	(1315,445)(1329,456)(1345,466)
	(1361,475)(1378,485)(1395,493)
	(1413,501)(1431,509)(1450,516)
	(1466,521)(1483,526)(1501,532)
	(1519,537)(1538,542)(1558,548)
	(1578,554)(1599,561)(1620,569)
	(1641,577)(1662,586)(1682,595)
	(1702,606)(1720,617)(1738,629)
	(1755,642)(1770,656)(1784,670)
	(1797,687)(1809,704)(1819,721)
	(1828,739)(1837,759)(1844,780)
	(1852,803)(1858,827)(1865,852)
	(1870,878)(1876,905)(1881,933)
	(1885,961)(1889,989)(1893,1017)
	(1896,1045)(1899,1072)(1902,1099)
	(1905,1124)(1908,1149)(1911,1172)
	(1913,1194)(1916,1215)(1919,1235)
	(1924,1261)(1929,1285)(1935,1308)
	(1942,1330)(1950,1352)(1959,1374)
	(1969,1397)(1981,1420)(1993,1442)
	(2005,1464)(2017,1483)(2026,1499)
	(2033,1510)(2037,1517)(2039,1520)
\put(209,27){\makebox(0,0)[lb]{\smash{{{\SetFigFont{8}{9.6}{\rmdefault}{\mddefault}{\updefault}0}}}}}
\put(82,2015){\makebox(0,0)[lb]{\smash{{{\SetFigFont{8}{9.6}{\rmdefault}{\mddefault}{\updefault}1}}}}}
\put(1986,50){\makebox(0,0)[lb]{\smash{{{\SetFigFont{8}{9.6}{\rmdefault}{\mddefault}{\updefault}1}}}}}
\put(0,1128){\makebox(0,0)[lb]{\smash{{{\SetFigFont{8}{9.6}{\rmdefault}{\mddefault}{\updefault}$\theta_2$}}}}}
\put(1028,50){\makebox(0,0)[lb]{\smash{{{\SetFigFont{8}{9.6}{\rmdefault}{\mddefault}{\updefault}$\theta_1$}}}}}
\put(5682,49){\makebox(0,0)[lb]{\smash{{{\SetFigFont{8}{9.6}{\rmdefault}{\mddefault}{\updefault}0}}}}}
\put(5532,274){\makebox(0,0)[lb]{\smash{{{\SetFigFont{8}{9.6}{\rmdefault}{\mddefault}{\updefault}0}}}}}
\put(5532,1999){\makebox(0,0)[lb]{\smash{{{\SetFigFont{8}{9.6}{\rmdefault}{\mddefault}{\updefault}1}}}}}
\put(7422,57){\makebox(0,0)[lb]{\smash{{{\SetFigFont{8}{9.6}{\rmdefault}{\mddefault}{\updefault}1}}}}}
\put(2977,42){\makebox(0,0)[lb]{\smash{{{\SetFigFont{8}{9.6}{\rmdefault}{\mddefault}{\updefault}0}}}}}
\put(4702,42){\makebox(0,0)[lb]{\smash{{{\SetFigFont{8}{9.6}{\rmdefault}{\mddefault}{\updefault}1}}}}}
\put(2819,252){\makebox(0,0)[lb]{\smash{{{\SetFigFont{8}{9.6}{\rmdefault}{\mddefault}{\updefault}0}}}}}
\put(2827,1970){\makebox(0,0)[lb]{\smash{{{\SetFigFont{8}{9.6}{\rmdefault}{\mddefault}{\updefault}1}}}}}
\put(3734,50){\makebox(0,0)[lb]{\smash{{{\SetFigFont{8}{9.6}{\rmdefault}{\mddefault}{\updefault}$\tilde\theta_1$}}}}}
\put(6434,50){\makebox(0,0)[lb]{\smash{{{\SetFigFont{8}{9.6}{\rmdefault}{\mddefault}{\updefault}$\theta'_1$}}}}}
\put(5444,1113){\makebox(0,0)[lb]{\smash{{{\SetFigFont{8}{9.6}{\rmdefault}{\mddefault}{\updefault}$\theta'_2$}}}}}
\put(2744,1077){\makebox(0,0)[lb]{\smash{{{\SetFigFont{8}{9.6}{\rmdefault}{\mddefault}{\updefault}$\tilde\theta_2$}}}}}
\put(6792,58){\makebox(0,0)[lb]{\smash{{{\SetFigFont{8}{9.6}{\rmdefault}{\mddefault}{\updefault}$\frac1{\gamma}$}}}}}
\end{picture}
}
\end{center}
\caption{The action of $T$ and $U$ on an orbit of the lifted flow of $X$.}     
\label{action of T}
\end{figure}

Considering a perturbation of a constant ``irrational'' vector field
$\omega$ lifted to $\Rr^2$, we
want to determine under which conditions there is still equivalence to
$\omega$.
Let $X$ be a vector field close to $\omega$ arising from a time-independent
perturbation. 
We choose e.g. $\omega=(1, \gamma)$ where $\gamma=\frac{1+\sqrt {5}}{2}$
is the golden ratio.
The main idea is to perform a change of basis, from the canonical base to
$\{(0,1),(1,1)\}$, enlarging the region of $\Rr^2$ around the orbits
of $\omega$
(see Figure \ref{action of T}).
This is achieved by the 
linear transformation:
$\theta'=T^{-1}\theta$,
where 
$T=\left[\begin{smallmatrix}0&1\\1&1\end{smallmatrix}\right]\in GL(2,\Zz)$. 
The eigenvalues and eigenvectors of $T$ are given by:
$T\omega=\gamma\omega$ and $T\Omega=-\frac1\gamma\Omega$, where $\Omega\perp\omega$.
We also rescale the time, $t'=\frac1\gamma t$, because the orbits
take longer to cross the new torus.
This transformation does not affect the unperturbed vector field $\omega$,
as it is given by 
$\gamma DT^{-1}\circ T(\theta')\,\omega\circ T (\theta')=\gamma
T^{-1}\omega=\omega$.

We should now consider a coordinate
change $h(\theta')=U_X\circ T(\theta')$ with $U$ satisfying $U_{\omega}=\id$.
The fundamental requirement is that it has to cancel the growth of the
wiggles in the orbits of $X$, which are enlarged by
the use of $T$
(see orbit in Figure \ref{action of T}).
The idea is to perform a non-linear,
close to the identity periodic coordinate transformation $U_X$, in order to
remove as many perturbation terms as possible.

There are some terms which are more relevant -- ``resonant''.
To understand what they are and how they appear, consider the vector field $X$ in the
form of its Fourier decomposition:
$$
\dot \theta = \omega + \varepsilon \sum\limits_{k\in \Zz^2}f_k
e^{2\pi ik\cdot\theta }
$$
with $f_k\in\Cc^2$ and $\varepsilon>0$ ``small''.
We can approximately solve the above system using the unperturbed
solution plus an order $\varepsilon$ term. So, on the universal cover,
\begin{equation}\label{perturbed solution}
\theta(t)=\theta_0+\omega t + \varepsilon f_0 \,t +\varepsilon 
\sum\limits_{k\in \Zz^2\setminus\{0\}}  \frac{f_k}{2\pi ik\cdot\omega} e^{2\pi
ik\cdot(\theta_0+\omega t)} + \OO(\varepsilon^2), 
\end{equation}
for $t>0$.
The resonant terms are the ones whose index $k$ is
almost perpendicular to $\omega$.

We need to add some conditions to $\omega$ in order to avoid having small
denominators in the solution of the flow. 
This is done in the usual KAM-type proofs
by imposing a diophantine condition on $\omega$, i.e. 
the denominators $\abs{k\cdot\omega}$ in
(\ref{perturbed solution}) admit a
lower bound.
We will show in the following that it is enough to find $U$ such that
it eliminates only ``far from resonance terms''.
This is so because of the extra linear change of coordinates $T$
described before, that is responsible for a ``shift'' of resonant into
non-resonant modes.
In some cases, by iterating this process the orbits can be straightened.

We define $U$ for any non-zero vector $\omega$.
Naturally, a restriction on $\omega$ will appear again in the procedure
with $T$ (see also \cite{jld2}), since this type of results fails for some vectors with irrational
slope, in particular for some non-diophantine vectors
\cite{Arnold2}.

%%%%%%%%%%%%%%%%%%%%%%%%%%%%%%%%%%%%%%%%%%%%%%%%%%%%%%%%%%%%%%%%%%%%%%%%%%%%

\section{Space of Analytic Vector Fields}\label{Space of Analytic Vector Fields}

The following is valid for any $d$-dimensional torus $\Tt^d$,
$d\geq2$, for a class of frequencies $\omega$ to be defined.
The vector fields considered are inside a ball around $\omega$ in some
adequate space, and can
be regarded as maps of $\Rr^d$ by lifting their domains.
We make use of the analyticity to extend to the complex domain, so we
deal with complex analytic vector fields. 
We construct the
renormalisation operator $\RR$ and we look at the spectral properties
of its derivative at the fixed point $\omega$, to relate to the local dynamics.

The use of analytic function spaces as the domain of the renormalisation
functional operator $\RR$ is justified by the
usefulness of $\RR$ being $C^1$. 
Otherwise the ``picture'' of $\RR$ being a dynamical
system with stable and unstable manifolds of the fixed points would vanish.
The problem lies in the compositions that appear in $\RR$, thus
reducing the degree of differentiability of the image of the
renormalisation operator.
Considering $C^{\infty}$ functions would also be a possibility, but
complicate the technical parts of the method.

\medskip

Let $r>0$ and the domain lifted to a complex
neighbourhood of $\Rr^d$:
$$
\DD(r)=\left\{\theta\in\Cc^d \colon \|\im\theta\| <\frac r{2\pi}\right\},
$$
where $\|\cdot\|$ is the $\ell_1$-norm on $\Cc^d$.
That is,
$\|z\|=\sum_{i=1}^{d}|z_i|$, with 
$|\cdot|$ the usual norm on $\Cc$.
We will also be using the inner product:
$z\cdot z'=\sum_{i=1}^{d}z_iz_i'$.

An analytic function $f\colon\DD(r)\to\Cc^d$, $2\pi$-periodic in each
variable $\theta_i$, is represented in
Fourier series as
$$
f(\theta)=\sum\limits_{k\in\Zz^d}f_ke^{2\pi i k\cdot \theta},
$$
with coefficients $f_k\in\Cc^d$.
The Banach spaces $(\A_d(r),\|\cdot\|_r)$ and
$(\A'_d(r),\|\cdot\|'_r)$ are the subsets of the set of these
functions
such that the respective norms 
$$
\|f\|_r=\sum\limits_{k\in\Zz^d}\|f_k\|e^{r\|k\|}
\quad
\text{ and }
\quad
\|f\|'_r=\sum\limits_{k\in\Zz^d} \left(1+2\pi\|k\|\right)\|f_k\| e^{r\|k\|}
$$
are finite.

Let $\rho>0$ be a fixed value in the following.
Consider the vector fields of the form
$X(\cdot)=\omega+f(\cdot)$, where $\omega\in\Rr^d$ 
and $f\in \A_d'(\rho)$.
Assume that $X$ has
no equilibrium points, $X(\theta)\not=0$ for any
$\theta\in\Tt^d$.
A condition like $\|f\|'_\rho<\|\omega\|$
is enough to assure that there are no equilibria and, for any $\theta\in\Tt^d$, there is no
$\theta'\in\Tt^d$ such that $X(\theta')/\|X(\theta')\|=-X(\theta)/\|X(\theta)\|$. 
So, we could not have two symmetric winding ratios
for
such a vector field, since it would have implied having at
least two symmetric normalised tangent vectors.    
The set of such vector fields on $\Tt^2$ corresponds to the class of
the Poincar\'e flows\index{Poincar\'e flows} 
(fixed-point-free with only one possible winding ratio, see e.g. \cite{BGKM}).
We study these systems in more detail in \cite{jld2}.

%%%%%%%%%%%%%%%%%%%%%%%%%%%%%%%%%%%%%%%%%%%%%%%%%%%%%%%%%%%%%%%%%%%%%%%%%%%%

\section{Elimination of the Far from Resonance Terms}\label{section
elimination far from resonance terms for local Td}

We will be dealing with the vector fields considered above decomposed
in Fourier series:
$X(\theta)=\sum_{k\in \Zz^d}X_ke^{2\pi ik\cdot\theta}$.
As in \cite{Koch}, it is important to distinguish different
classes of terms according to their respective indices.

\medskip

\begin{definition}\label{far from resonance terms}
For a fixed value of $\sigma>0$ and
$\psi\in\Cc^d$, we define the {\it far
from resonance terms}\index{far
from resonance terms} with respect to $\psi$ to be the ones whose indices are in
$$
I_\sigma^-(\psi)=\left\{ k\in\Zz^d\colon |\psi\cdot k| > 
\sigma\|k\|\right\}.
$$
Similarly, the {\it resonant terms}\index{resonant terms} are $I_\sigma^+(\psi)=\Zz^d\setminus I_\sigma^-(\psi)$.
It is also useful to define the projections
$\Ii_\sigma^+(\psi)$ and $\Ii_\sigma^-(\psi)$
for any vector field $X$ by
$$
[\Ii_\sigma^\pm(\psi)] X(\theta)=\sum\limits_{k\in I_\sigma^\pm(\psi)}X_ke^{2\pi ik\cdot\theta}
$$ 
and $\Ii=\Ii_\sigma^+(\psi)+\Ii_\sigma^-(\psi)$ is the identity operator.
\end{definition}

The existence of a nonlinear change of coordinates $U$, close to the
identity, that eliminates the far from resonance terms $I^-_\sigma(\omega)$ of
a vector field $X$ in a neighbourhood of a non-zero constant vector field, is given by

\begin{theorem}\label{main theorem1}
Let $0<\rho'<\rho$, $\omega\in \Cc^d\setminus\{0\}$ and
$0<\sigma<\|\omega\|$.
If $X$ is in the open ball $\hat B\subset\A'_d(\rho)$ centred at
$\omega$ with radius 
$$
\hat\varepsilon=
\frac{\sqrt6-2}{12}\sigma
\min\left\{
\frac{\rho-\rho'}{4\pi},
\frac{3-\sqrt6}{6} \frac\sigma{\|\omega\|}
\right\},
$$ 
there is an invertible coordinate change $U\colon\DD(\rho')\to\DD(\rho)$
close to the identity, satisfying 
$$
\begin{array}{c}
[\Ii_\sigma^-(\omega)](DU)^{-1}\,X\circ U=0,
\quad\text{ and }\quad
U=\id 
\text{ if }[\Ii_\sigma^-(\omega)]X=0.
\end{array}
$$ 
The map $\UU\colon \hat B \to [\Ii_\sigma^+(\omega)]\A_d(\rho')$ given by $X\mapsto (DU)^{-1}\,X\circ U$
is analytic, and the derivative at every constant vector field $\omega+\psi$
in $\hat B$ is equal to $\Ii_\sigma^+(\omega)$.
Moreover, 
$$
\|\UU(X)-\omega\|_{\rho'}\leq
2\left(
1+\max\left\{
\frac23(3-\sqrt6),6 (\sqrt6+2)\frac{\|\omega\|}\sigma
\right\}
\right)
\|X-\omega\|'_\rho.
$$ 
\end{theorem}

In Section \ref{section proof of
main theorem} we include a proof of the above theorem using the
``homotopy method''.
It is worth noting here that since we are requiring to eliminate only the
far from resonance terms outside $I^+_\sigma(\omega)$, we avoid the
problem in the proof of dealing with small
denominators\index{small
denominators} (see Section \ref{proof of prop}).

%%%%%%%%%%%%%%%%%%%%%%%%%%%%%%%%%%%%%%%%%%%%%%%%%%%%%%%%%%%%%%%%%%%%%%%%%%%%

\section{Change of Basis}\label{Change of Basis}

All the coordinate changes of $\Tt^d$ are isotopic to a linear transformation in
$GL(d,\Zz)$ which preserves the structure and volume of the $d$-torus.
That can be seen as a map on the frequency vectors space, acting
on a vector field by shifting terms in the space of indices $\Zz^d$
and changing coefficients.

\medskip

\begin{definition}
For a fixed non-zero vector $\omega\in\Rr^d$
with rationally
independent components, i.e. $k\cdot\omega\not=0$ for any non-zero integer
vector $k$,
assume that $\omega=(1,\omega_2,\dots,\omega_d)$ where its
components are real algebraic numbers (roots of non-zero polynomials
over $\Qq$). 
The vectors $c\,\omega$, $c\in\Rr\setminus\{0\}$, are said to be of
{\it Koch type}, i.e. to belong to the set $KT_d$ if
the algebraic extension of $\Qq$ by the numbers
$\omega_2,\dots,\omega_d$ (the smallest field containing $\omega_2,\dots,\omega_d$ and
$\Qq$) is of degree $d$.
\end{definition}

An important property, 
the existence of a linear change of basis $T$ corresponding to the above
class of vectors, follows from

\begin{lemma}[Koch \cite{Koch} -- Lemma 4.1]\label{exist T}
A vector $\omega$ is in $KT_d$ 
if and only if
there is an integral $d\times d$ matrix $T$, such that
$0<|\lambda_d|\leq \dots \leq |\lambda_2|<1<|\lambda_1|$ and 
$T\omega=\lambda_1\omega$,
where $\lambda_1, \dots, \lambda_d$ are simple eigenvalues of $T$.
Furthermore, some of the (integral) matrices satisfying these
properties are in $GL(d,\Zz)$.
\end{lemma}

\begin{remark} \label{remark algebraic extension}
For $d=2$, the numbers $\omega_2$ that produce a degree 2 algebraic
extension of $\Qq$ are the quadratic irrationals\index{quadratic irrational} (roots of
non-zero degree 2 polynomials over $\Qq$). These are characterised by an
eventually periodic continued fraction expansion \cite{Lang}.
For $d=3$, the base for the algebraic extension of $\Qq$ by $\alpha$
and $\beta$ is $\{1,\alpha,\alpha^2\}$, where $\alpha$ has to be a
cubic irrational\index{cubic irrational} (root of a non-zero degree 3 polynomial over $\Qq$).
Then $\beta=c_1+c_2\alpha+c_3\alpha^2$, $c_i\in\Qq$.
\end{remark}

In what follows, a vector $\omega\in KT_d$ and a corresponding matrix $T\in
GL(d,\Zz)$ are chosen according to Lemma \ref{exist T}. 
Also, $\bar\omega\in\Rr^d$ is chosen to be the (unstable)
$\lambda_1$-eigenvector of $T$ such that 
$\bar\omega\cdot\omega=1$.
Notice that $KT_d$ is a subset of the diophantine vectors
$DC(\beta)$ with $\beta=-1-\ln|\lambda_1|/\ln|\lambda_2|$ (\cite{Koch}
-- Corollary 4.2),
i.e. there is a
constant $C>0$ that satisfies
$|\omega\cdot k| > C \|k\|^{-\beta-1}$ for every integer vector $k\not=0$.
In addition,
$\prod_{j=1}^d |\lambda_j|=1$,
where $\lambda_j$ are the eigenvalues of $T$, with $\lambda_1\in\Rr$ since it is the only one outside the unit
circle.
Also, $\omega\cdot \omega^{(j)}=0$, with $\omega^{(j)}$ being
the eigenvector of $T$ corresponding to $\lambda_j$, $2\leq j\leq d$,
because $\omega\cdot\omega^{(j)}=T^*\omega\cdot
T^{-1}\omega^{(j)}=(\lambda_1/\lambda_j)\omega\cdot\omega^{(j)}$.
We denote the transpose matrix of $T$ by $T^*$.

From now on, we simplify the notation
writing $I^\pm_\sigma$ and $\Ii^\pm_\sigma$ instead of the
sets $I^\pm_\sigma(\omega)$ and the projections
$\Ii^\pm_\sigma(\omega)$, respectively.

The next proposition determines that the change of basis corresponding to the matrix $T$ above, acting on the space
of ``resonant vector fields'', is analyticity improving\index{analyticity improving}.

\begin{proposition}\label{proposition T}
Let $0<\kappa<1$ and $0<\rho'<\rho$. 
If $\kappa\rho<\rho'$, then,
for some $0<\sigma<\frac12\|\bar\omega\|^{-1}$, any $X\in\Ii_\sigma^+\A_d(\rho')$
has an analytic extension on
$T\DD(\rho)$. The linear map $\TT(X)= T^{-1}\,X\circ T$ from $\Ii_\sigma^+\A_d(\rho')$
to $\A'_d(\rho)$ is compact. 
\end{proposition}

The proof is given in Section \ref{Analyticity Improvement}.
The values of $\sigma$ have to be chosen
such that
$$
I^+_\sigma\subset I^\kappa:=\{x\in\Rr^d\colon\|T^*x\|\leq\kappa\|x\|\}
$$ 
for a
given $\kappa<1$.
The proof of this proposition only requires that
the resonant terms are inside the set $I^\kappa$.
Thus, we could have different definitions for $I^+_\cdot$.
The different choices are restricted
by the part of the proof of Theorem \ref{main theorem1} in Section
\ref{proof of prop}, that requires the
denominator terms $|\omega\cdot k|$ to have at least a $\|k\|$-linear
lower bound on $I^-_\cdot$.

It is also important to remark that one assumes that the orthogonal
hypersurface in $\Rr^d$ with respect to $\omega$ is contained
in $I^\kappa$. 
%That is, $\max \{\|T^*y\|\colon \|y\|=1,\omega\cdot y=0\}<\kappa$.
This condition is indeed verified by choosing appropriate norms on $\Rr^d$
or by considering a sufficiently large power of $T$ instead.
The indices of resonant terms are therefore shifted
towards $I^-_\sigma$, eventually becoming far from resonance terms.

%%%%%%%%%%%%%%%%%%%%%%%%%%%%%%%%%%%%%%%%%%%%%%%%%%%%%%%%%%%%%%%%%%%%%%%%%%%%

\section{Renormalisation Operator}\label{section renorm oper}

We are now in condition to construct a renormalisation operator $\RR$
based on the coordinate transformations introduced before.

\medskip

Fix $\omega\in KT_d$, $T$ and $\lambda_j$, $1\leq j\leq d$, as in Lemma \ref{exist T}.
We denote by $\Ee(X)$ the average on $\Tt^d$ of a vector field
$X$, i.e.
$$
\Ee(X)=\int_{\Tt^d}X(\theta)\,d\theta,
$$
where $d\theta$ is the
normalised Lebesgue measure on $\Tt^d$.

\begin{definition}
Consider the maps $\UU$ and $\TT$ given by Theorem \ref{main theorem1}
and Proposition \ref{proposition T}, respectively.
The {\it renormalisation operator}\index{renormalisation operator} $\RR$ is
defined to be
$$
\RR(X) = \frac{1}{\bar\omega\cdot\Ee\,(\tilde X)} \tilde X,
\quad \tilde X = \TT\circ\UU(X), 
$$
for vector fields $X$.
\end{definition}

The rescaling of time is chosen to preserve the term $\omega$, because
without it the iterates of $\omega$, 
$(\TT\circ\UU)^n(\omega)=\lambda^{-n}_1\omega$, $n\geq0$, would go to zero.
We could fix that by simply considering the product by $\lambda_1$ inside the
calculation of $\tilde X$, but then
all the elements of the one-parameter family $X_\mu=(1+\mu)\,\omega$,
$\mu\in\Cc$,
would be fixed points of the renormalisation.
This ``neutral'' direction is contracted by the operator $\RR$ given in the definition, 
allowing $\omega$ to be an isolated fixed point.

This choice of rescaling has the disadvantage of reducing the
domain of the operator (see Proposition \ref{R well-defined} below),
when comparing with the size of the ball $\hat B$ obtained
in Theorem \ref{main theorem1}.
In any case, this can be solved by using the rescaling suggested above, thus
working with a one-parameter family of fixed
points and a larger domain.

Note also that, by Proposition \ref{proposition T},
$\sigma<\frac12 \|\omega\|^{-1} \|\omega\|^2_2\leq\frac12 \|\omega\|$, where
$\|\cdot\|_2$ represents the Euclidean norm.
This means that the radius of an open ball given by
Theorem \ref{main theorem1} can be taken to be $\CC\sigma^2\|\omega\|^{-1}$, assuming
at least $\CC< 4$.
Hence, the domain of $\RR$ only contains equilibria-free vector fields.

\begin{figure}
\begin{center}
\input{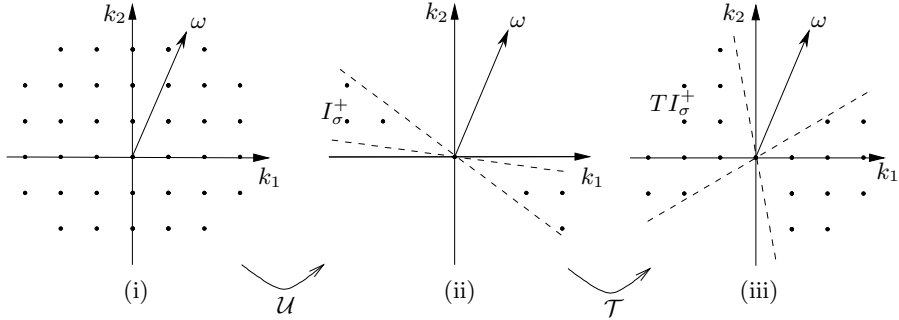} 
\end{center}
\caption{Steps of the renormalisation operator $\RR$ for $d=2$, on the
lattice $\Zz^2$, the indices of the
Fourier coefficients (notice that each coefficient is a vector
in $\Cc^d$).}     
\label{renormalisation steps}
\end{figure}

The change of basis given by $T$ turns some
of the resonant terms into $I_\sigma^-$
where they will be completely eliminated by the operator $\UU$ of the
next iteration of $\RR$.
The use of such $\UU$ is not essential, it would be enough to
``sufficiently'' reduce
the resonant terms in a way that analyticity
would not be lost.
As the elimination is complete, we get a faster convergence to the
fixed point and a much simpler and clearer analysis starting from the
fact that the derivative of $\UU$ is straightforward.
If we would only apply $\UU$, all the neutral directions that are
eliminated by the use of $\TT$ would persist (no elimination of the
resonant terms), and the analyticity domain would
not be fully recovered.
The right combination of the two coordinate changes is
responsible for the usefulness of $\RR$ (see Figure
\ref{renormalisation steps} for the $d=2$ case).

\begin{proposition}\label{R well-defined}
Let $\rho>0$. 
The renormalisation operator $\RR$ is a well-defined analytic map from
an open neighbourhood $B\subset\A_d'(\rho)$ of $\omega$, to
$\A'_d(\rho)$.
\end{proposition}

\proof
Consider the complex-valued continuous functional $F(\tilde
X)=\bar\omega\cdot\Ee(\tilde X)$,
and restrict its domain to a neighbourhood
$\tilde B\subset\A_d'(\rho)$ of $\lambda_1^{-1}\omega$ such that
$F(\tilde B)\subset\{z\in\Cc\colon|\lambda_1 z-1|<1/2\}$, i.e. 
$F$ is bounded away from zero in $\tilde B$.
 
The domain $B$ of $\RR$ is a subset of $\hat B$ (the domain of $\UU$) such that $\TT\circ
\UU(B) \subset \tilde B$.
That is satisfied for a small enough choice of the radius of
$B\subset\hat B$,
since
$\UU(X)=\omega+\Ii^+_\sigma f+\OO({\|f\|'_\rho}^2)$ and
$|\lambda_1 F(\TT\circ\UU(X))-1| = 
|\lambda_1 (\bar \omega \cdot T^{-1}\,\Ee\,f)+ \OO({\|f\|'_\rho}^2)|$,
with $X=\omega+f$.
It is then sufficient to have $|\bar\omega \cdot \Ee\, f +
\OO({\|f\|'_\rho}^2) | <\frac12$.
Notice that $\bar\omega\cdot T^{-1}\Ee\,f=T^{-1*}\bar\omega\cdot\Ee\,f=\lambda_1^{-1}(\bar\omega\cdot\Ee\,f)$.
Hence, as
$F$ is bounded and analytic in
$\tilde B$,
the above, Theorem \ref{main theorem1} and Proposition
\ref{proposition T} prove the claim.
\cqd

%%%%%%%%%%%%%%%%%%%%%%%%%%%%%%%%%%%%%%%%%%%%%%%%%%%%%%%%%%%%%%%%%%%%%%%%%%%%

\section{Hyperbolicity of the Fixed Point $\omega$}\label{Hyperbolicity of the Fixed Point omega}

The transformation $\RR$ was constructed in a way to hold an isolated fixed
point at $\omega$.
The local behaviour of $\RR$ around it is given by the derivative,
with which we hope to characterise the vector fields close
enough to $\omega$, realising the existence or not of an
analytic equivalence.

\medskip

We can rewrite $\RR$ by
the expression $\RR=\FF\circ\TT\circ\UU$,
where the time-rescaling step is given by 
$\FF(X)=[\bar\omega\cdot\Ee\,(X)]^{-1}X$ 
and $D\FF(\lambda_1^{-1}\omega)f=\lambda_1^{-1}[f-\bar\omega\cdot\Ee\,(f)]\,\omega$.
The derivative of $\RR$ at $\omega$ is then the linear map
$$
D\RR(\omega)\, f=  \LL\,f-[\bar\omega\cdot\Ee\,(\LL\,f)]\,\omega, 
\quad 
f\in\A'_d(\rho),
\quad \text{
with }\quad
\LL\,f=\lambda_1 \TT \circ\Ii_\sigma^+ f.
$$ 
This operator is compact by Proposition \ref{proposition T} and using
the
fact that $\Ii^+_\sigma$ is bounded.
The eigenvalues of $\LL$ are zero and those of $T^{-1}$ multiplied by $\lambda_1$.
The vector $\omega$ gives the subspace
corresponding
to the eigenvalue $1$.
Since the renormalisation
should not depend on the size of $\omega$, this direction is fully 
eliminated in $D\RR(\omega)$ by the rescaling of
time: $D\RR(\omega)\,\omega=\omega-(\bar\omega\cdot\omega) \, \omega=0$.
There are $(d-1)$ unstable directions of $D\RR(\omega)$ given by the other eigenvectors
of $T^{-1}$,
with eigenvalues $\lambda_1\lambda_j^{-1}$,
$2\leq j \leq d$, of modulus greater than one
(see Section \ref{Spectral Properties and Invariant Manifolds}).

The remaining eigenvalue of $D\RR(\omega)$ is zero, meaning that the
elimination of only the far from resonant terms in conjuction with the
change of basis are sufficient to eliminate
all the non-constant Fourier terms of a vector field close to $\omega$.
Therefore, there is a codimension-$(d-1)$ manifold inside a neighbourhood of
$\omega$ in $\A'_d(\rho)$ being mapped by $\RR$ into itself, and
we have proved

\begin{theorem}\label{renomalisation picture th1}
If $\omega\in KT_d$, then the constant vector field $\omega$ is a hyperbolic
fixed point of $\RR$ with a local codimension-$(d-1)$
stable manifold $\WW^s(\omega)$ and a local $(d-1)$-dimensional unstable manifold $\WW^u(\omega)$. 
\end{theorem}

\begin{figure}
\begin{center}
\setlength{\unitlength}{0.00062500in}
\begingroup\makeatletter\ifx\SetFigFont\undefined%
\gdef\SetFigFont#1#2#3#4#5{%
  \reset@font\fontsize{#1}{#2pt}%
  \fontfamily{#3}\fontseries{#4}\fontshape{#5}%
  \selectfont}%
\fi\endgroup%
{\renewcommand{\dashlinestretch}{30}
\begin{picture}(4223,2328)(0,-10)
\put(1291,1139){\blacken\ellipse{20}{20}}
\put(1291,1139){\ellipse{20}{20}}
\path(935,1195)(1000,1230)(955,1275)
\path(1722,1291)(1662,1311)(1687,1381)
\path(1741,970)(1691,1045)(1766,1065)
\path(1291,2020)(1292,2049)
\path(1317.847,1928.037)(1292.000,2049.000)(1257.882,1930.105)
\path(1301,357)(1301,327)
\path(1271.000,447.000)(1301.000,327.000)(1331.000,447.000)
\path(936,1032)(998,1082)(936,1112)
\thicklines
\path(1257,1551)(1255,1551)(1249,1551)
	(1239,1550)(1224,1549)(1203,1548)
	(1176,1547)(1144,1545)(1108,1543)
	(1067,1540)(1025,1538)(980,1535)
	(936,1532)(892,1529)(849,1525)
	(807,1522)(768,1519)(731,1515)
	(696,1512)(664,1508)(633,1504)
	(605,1500)(578,1496)(552,1492)
	(528,1488)(505,1483)(482,1478)
	(460,1473)(433,1465)(407,1458)
	(381,1449)(355,1440)(330,1430)
	(305,1420)(280,1408)(256,1396)
	(233,1384)(211,1370)(189,1356)
	(169,1342)(150,1327)(132,1311)
	(115,1296)(100,1280)(86,1264)
	(74,1248)(63,1232)(54,1216)
	(46,1200)(39,1183)(34,1171)
	(30,1158)(27,1145)(24,1132)
	(23,1119)(22,1105)(23,1091)
	(24,1076)(27,1061)(31,1046)
	(36,1031)(43,1016)(52,1000)
	(61,985)(73,969)(86,953)
	(101,938)(117,923)(135,908)
	(155,893)(177,879)(200,865)
	(224,852)(251,839)(278,826)
	(308,814)(339,802)(372,791)
	(407,780)(444,770)(472,762)
	(503,755)(534,747)(567,740)
	(601,733)(637,726)(674,720)
	(712,713)(752,707)(792,701)
	(834,695)(878,690)(922,684)
	(967,679)(1013,675)(1060,670)
	(1107,666)(1155,662)(1203,659)
	(1252,656)(1300,653)(1349,651)
	(1397,649)(1445,647)(1492,646)
	(1539,646)(1585,646)(1630,646)
	(1675,647)(1718,648)(1761,650)
	(1802,652)(1842,654)(1881,657)
	(1918,660)(1955,664)(1990,668)
	(2024,673)(2057,678)(2089,683)
	(2129,691)(2168,700)(2205,710)
	(2241,721)(2275,732)(2308,745)
	(2340,758)(2370,772)(2399,787)
	(2427,803)(2453,819)(2477,836)
	(2501,854)(2522,872)(2543,890)
	(2561,909)(2578,928)(2594,947)
	(2608,966)(2620,985)(2631,1003)
	(2641,1022)(2650,1040)(2657,1057)
	(2663,1075)(2668,1092)(2672,1108)
	(2675,1124)(2677,1139)(2679,1154)
	(2679,1173)(2679,1191)(2678,1208)
	(2675,1226)(2672,1243)(2666,1259)
	(2660,1276)(2652,1292)(2642,1308)
	(2631,1323)(2619,1338)(2605,1352)
	(2590,1366)(2573,1379)(2555,1392)
	(2536,1403)(2516,1414)(2494,1425)
	(2472,1434)(2448,1443)(2424,1451)
	(2399,1459)(2372,1466)(2345,1473)
	(2323,1477)(2301,1482)(2277,1486)
	(2252,1490)(2226,1494)(2198,1497)
	(2169,1501)(2137,1504)(2103,1508)
	(2067,1511)(2029,1514)(1988,1517)
	(1944,1520)(1899,1523)(1851,1526)
	(1801,1529)(1750,1531)(1698,1534)
	(1646,1537)(1596,1539)(1547,1541)
	(1501,1543)(1459,1545)(1421,1547)
	(1390,1548)(1364,1549)(1344,1550)
	(1329,1550)(1320,1551)(1315,1551)(1313,1551)
\thinlines
\path(1278,1136)(1275,1137)(1267,1140)
	(1255,1145)(1237,1152)(1216,1160)
	(1193,1169)(1169,1178)(1146,1186)
	(1124,1194)(1103,1200)(1083,1206)
	(1064,1212)(1045,1216)(1026,1220)
	(1007,1224)(991,1227)(975,1230)
	(957,1233)(938,1235)(917,1238)
	(895,1240)(870,1243)(843,1245)
	(814,1248)(783,1251)(751,1254)
	(718,1256)(687,1259)(659,1261)
	(635,1263)(616,1264)(604,1265)
	(596,1266)(593,1266)
\path(2030,876)(2027,878)(2022,881)
	(2011,887)(1997,896)(1979,907)
	(1958,920)(1935,933)(1912,946)
	(1889,959)(1868,970)(1847,981)
	(1827,991)(1809,1000)(1790,1008)
	(1772,1016)(1754,1023)(1735,1029)
	(1716,1036)(1695,1042)(1674,1049)
	(1652,1055)(1629,1061)(1606,1067)
	(1582,1073)(1558,1079)(1534,1085)
	(1510,1090)(1487,1095)(1465,1100)
	(1444,1104)(1425,1108)(1407,1112)
	(1391,1115)(1377,1118)(1363,1120)
	(1339,1125)(1321,1129)(1307,1131)
	(1297,1133)(1290,1135)(1287,1136)(1285,1136)
\path(1296,1147)(1299,1148)(1305,1152)
	(1316,1157)(1331,1165)(1352,1175)
	(1376,1187)(1402,1200)(1430,1214)
	(1457,1228)(1484,1240)(1510,1252)
	(1534,1263)(1556,1273)(1576,1282)
	(1596,1290)(1614,1297)(1631,1303)
	(1648,1309)(1665,1314)(1684,1320)
	(1703,1325)(1722,1329)(1743,1334)
	(1764,1338)(1788,1341)(1813,1345)
	(1841,1349)(1870,1352)(1900,1356)
	(1930,1359)(1958,1362)(1982,1364)
	(2002,1366)(2015,1367)(2023,1368)(2027,1368)
\path(1301,647)(1301,12)
\path(1291,2301)(1291,2298)(1291,2293)
	(1291,2282)(1291,2266)(1291,2244)
	(1291,2215)(1291,2181)(1291,2142)
	(1291,2099)(1291,2054)(1291,2008)
	(1291,1961)(1291,1914)(1291,1869)
	(1291,1827)(1291,1786)(1291,1748)
	(1291,1712)(1291,1679)(1291,1649)
	(1291,1620)(1291,1594)(1291,1569)
	(1291,1546)(1291,1525)(1291,1505)
	(1291,1485)(1291,1457)(1291,1430)
	(1291,1405)(1291,1380)(1291,1356)
	(1291,1333)(1291,1309)(1291,1285)
	(1291,1260)(1291,1236)(1291,1213)
	(1291,1192)(1291,1174)(1291,1160)
	(1291,1150)(1291,1145)(1291,1142)
\path(552,900)(555,901)(561,905)
	(572,910)(587,918)(608,928)
	(632,940)(658,953)(686,967)
	(713,981)(740,993)(766,1005)
	(790,1016)(812,1026)(832,1035)
	(852,1043)(870,1050)(887,1056)
	(904,1062)(921,1067)(940,1073)
	(959,1078)(978,1082)(999,1087)
	(1020,1091)(1044,1094)(1069,1098)
	(1097,1102)(1126,1105)(1156,1109)
	(1186,1112)(1214,1115)(1238,1117)
	(1258,1119)(1271,1120)(1279,1121)(1283,1121)
\put(1430,2012){\makebox(0,0)[lb]{\smash{{{\SetFigFont{9}{10.8}{\rmdefault}{\mddefault}{\updefault}$\WW^u(\omega)$}}}}}
\put(2731,1022){\makebox(0,0)[lb]{\smash{{{\SetFigFont{9}{10.8}{\rmdefault}{\mddefault}{\updefault}$\WW^s(\omega)$}}}}}
\put(1119,925){\makebox(0,0)[lb]{\smash{{{\SetFigFont{9}{10.8}{\rmdefault}{\mddefault}{\updefault}$\omega$}}}}}
\end{picture}
} 
\end{center}
\caption{The schematic renormalisation picture given by the action of $\RR$ in a
neighbourhood of $\omega$. Note that $\dim\WW^u(\omega)=d-1$ whilst
$\WW^s(\omega)$ is an infinite dimensional manifold (of codimension $d-1$).}     
\label{renormalisation picture}
\end{figure}
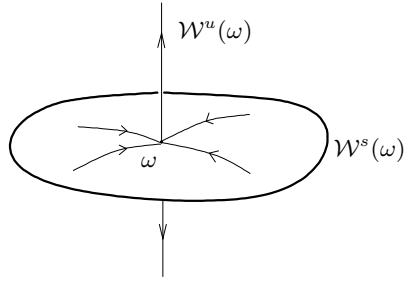

We have determined an analytic equivalence\index{analytic equivalence} between the elements of
$\WW^s(\omega)$ and in particular with $\omega$, that is given by
the operator $\lim\limits_{n\to+\infty}\RR^n$.
The local stable manifold $\WW^s(\omega)$ is the set of all $X$ in
some neighbourhood $B$ with
winding set $w_X=\omega/\|\omega\|$ (see Section \ref{Spectral Properties and Invariant Manifolds}).
On the other hand, the local unstable manifold $\WW^u(\omega)$ is the
affine space of the
vector fields in the form $X=\omega+v\in B$, where
$v\in\Span\{\omega^{(2)},\dots,\omega^{(d)}\}$, as this spectral
subspace is $\RR$-invariant.
A schematic representation of the renormalisation scheme described is
in Figure \ref{renormalisation picture}.

%%%%%%%%%%%%%%%%%%%%%%%%%%%%%%%%%%%%%%%%%%%%%%%%%%%%%%%%%%%%%%%%%%%%%%

\section{Proofs of Theorem \ref{main theorem1} and Proposition \ref
{proposition T}}\label{section proof of main theorem}
\label{Proofs of Theorems}

\subsection{Preliminaries}

We include here some technical details useful for the proofs.

Firstly, note that by writing $f=(f^1,\dots,f^d)$, one can rewrite the formulas of the
norms on the spaces $\A_d(r)$ and $\A_d'(r)$.
In fact, 
$$
\| f\|_r=\sum\limits_{i=1}^d\|f^i\|_r, \quad f^i\in\A_1(r)\quad
\text{ and }\quad
\| f\|'_r=\sum\limits_{i=1}^d\|f^i\|'_r,\quad f^i\in\A'_1(r),
$$
where
$$
\|f^i\|'_r = \|f^i\|_r + \sum\limits_{j=1}^d \|\partial_jf^i\|_r, \quad i=1,\dots,d.
$$
One has $\|f\|_r\leq\|f\|'_r$, $\A'_d(r)\subset\A_d(r)$ and
$\|f\|_{r'}\leq\|f\|_r$,  
for $0<r'<r$.

Define the inclusion maps $\II\colon\A_d(r)\to\A_d(r')$ and 
$\II'\colon\A'_d(r)\to\A'_d(r')$, $0<r'<r$, by the restriction
of the domain $\DD(r)$ to $\DD(r')$, i.e. $f\mapsto f|_{\DD(r')}$.
It is easy to check that $\II(\A_d(r)) \subset \A_d(r')$ and 
$\II'(\A'_d(r)) \subset \A'_d(r')$.
The inclusion maps $\II$ and $\II'$\index{inclusion map} are compact
linear operators.

Consider a bounded linear map $A\colon\A_d(r)\to\A_d(r)$ with
$(Af)(\theta)=A(\theta)\,f(\theta)$, where $A(\theta)$ is a linear map
on $\Cc^d$.
Its operator norm satisfies
$\|A\| \leq \sum_{i,j=1}^d \|a_{ij}\|_r$,
where $a_{ij}\in\A_1(r)$ such that, for each $\theta\in\DD(r)$, $a_{ij}(\theta)$ are the coefficients of
the matrix which represents $A(\theta)$ in the canonical basis.
In particular, if $A$ is the derivative of a function $f\in\A_d(r)$,
$(Df\,g)(\theta) = Df(\theta)\,g(\theta)$, $g\in\A_d(r)$ and
$\theta\in\DD(r)$, then 
$\|Df\|\leq\|f\|'_r$.

If $\alpha\in\A_1(r)$ and $f,g\in\A_d(r)$, then $\alpha\,
f\in\A_d(r)$, $f\cdot g\in\A_1(r)$,
$
\|\alpha\,f\|_r\leq\|\alpha\|_r\|f\|_r$ and $\|f\cdot g\|_r\leq\|f\|_r\|g\|_r.
$
Similarly, if $\alpha\in\A'_1(r)$ and $f,g\in\A'_d(r)$, then $\alpha\,
f\in\A'_d(r)$, $f\cdot g\in\A'_1(r)$,
$\|\alpha\,f\|'_r\leq\|\alpha\|'_r\|f\|'_r$ and $\|f\cdot g\|_r'\leq\|f\|'_r\|g\|'_r$.
In this way, the spaces $\A_d(r)$ and $\A_d'(r)$ are Banach algebras.

\begin{lemma}\label{lemma f-f*}
Let $0<r'<r$ and $f\in\A'_d(r)$.
If $u\in\A_d(r')$ and $\|u\|_{r'}<(r-r')/4\pi$,
then
\begin{enumerate}
\item\label{f*}
$\|f(\id+u)\|_{r'}\leq \|f\|_{(r+r')/2}$,
\item\label{Df*}
$\|Df(\id+u)\|\leq \|f\|'_{(r+r')/2}$,
\item\label{f*-f}
$\|f(\id+u)-f\|_{r'}\leq \|f\|'_{(r+r')/2}\,\|u\|_{r'}$,
\item\label{Df*-Df}
$\|Df(\id+u)-Df\|\leq \frac{4\pi}{r-r'}\|f\|'_{r}\,\|u\|_{r'}$.
\end{enumerate}
\end{lemma}

\proof
Knowing  that $\|e^{2\pi ik\cdot u}\|_{r'}\leq
e^{2\pi\|k\|\|u\|_{r'}}$ one gets \ref{f*} and \ref{Df*}.
The mean value theorem gives \ref{f*-f}.
To prove \ref{Df*-Df} apply again the mean value theorem now to
$\|\partial_jf^i(\id+u)-\partial_jf^i\|_{r'}$ and also
$$
\begin{array}{ll}
\|\partial_jf^i\|'_{(r+r')/2}&=\sum\limits_{k\in
\Zz^d}2\pi|k_j|(1+\|k\|)|f^i_k|e^{r\|k\|}e^{-(r-r')\|k\|/2}
\\
&\leq \frac{4\pi}{r-r'}\sum\limits_{k\in
\Zz^d}(1+\|k\|)|f^i_k|e^{r\|k\|}=\frac{4\pi}{r-r'}\|f^i\|'_r,
\end{array}
$$
where we have used the inequality: $\sup_{t\geq0}te^{-\beta
t}\leq1/\beta$ for $\beta>0$.
\cqd

%%%%%%%%%%%%%%%%%%%%%%%%%%%%%%%%%%%%%%%%%%%%%%%%%%%%%%%%%%%%%%%%%%%%%%%%%

\subsection{Homotopy Method\index{homotopy method}}\label{homotopy method}

The sometimes called ``homotopy trick''
has been used in different problems, such as by Moser to show that
all smooth volume forms on a compact orientable manifold are
equivalent up to a diffeomorphism 
(\cite{Katok} -- Section 5.1e).
Other examples of its application are
proofs of the Darboux theorem or the Poincar\'e lemma
(\cite{Katok} -- Sections
5.5.9 and A.3.11, respectively), and Roussarie's proof of Morse's lemma.
This procedure is also related to the ``deformation method''
used in different setups in KAM theory. 
It consists of a flow of symplectomorphisms that reduce in each
iterate the size of the perturbation of an integrable Hamiltonian
\cite{Poschel,Delshams}.

\medskip

In the following fix $\sigma>0$, thus dropping the index $\sigma$
for the sets $I^\pm$ and corresponding projections $\Ii^\pm$.
Choose also $\rho$ and $\rho'$ such that $0<\rho'<\rho$.
Let $0<\hat\delta<\frac12$ and $\hat\varepsilon$ be a
positive constant
which will be determined along the proof
and contains the restrictions on the size of the perturbation of $\omega$
depending on $\sigma$ and the norm of $\omega$, as will be seen in Section \ref{proof of prop}.

For vector fields in the form $X=\omega+f$, consider the open neighbourhood
$\EE$ of the term $f$:
$$
\EE=\{f\in\A'_d(\rho)\colon \|f\|'_{\rho}<  \hat\varepsilon\}.
$$
The coordinate transformation $U$ is written as $U=\id+u$, with $\id$
as the identity transformation, and $u$ in the open ball
$\BB$ of radius $\hat\delta$ in
$\Ii^-\A'_d(\rho')$, i.e.
$$
\BB=\left\{u\in\Ii^-\A'_d(\rho')\colon\|u\|'_{\rho'}
<\hat\delta\right\}.
$$

For a given function $f\in\EE$, consider the operator
$F\colon\BB\to\Ii^-\A_d(\rho')$,
$$
F(u)=\Ii^-(I+Du)^{-1}\,\left[\omega+f(\id+u)\right],
$$
where $I$ is the identity operator.
This is simply the transformed vector field $\Ii^-(DU)^{-1}\,X\circ U$, with $X=\omega+f$.

\begin{lemma}\label{DF formula}
The derivative of $F$ at $u\in\BB$ is the linear map from
$\Ii^-\A'_d(\rho')$ to $\Ii^-\A_d(\rho')$,
$$
\begin{array}{rl}
DF(u)\,h = & \Ii^-(I+Du)^{-1}\left[ Df\circ(\id+u)\,h 
\right.
\\
&
\left.
- Dh\,(I+Du)^{-1}\left(\omega+f\circ(\id+u)\right)
\right],
\end{array}
$$
with $h$ chosen such that $u+h\in\BB$.
\end{lemma}

\proof
We need to compute the linear term on $h$ of $F(u+h)-F(u)$.
As $u\in\BB$, we have that the bounded linear operator $Du$, from
$\A_d(\rho')$ into itself, satisfies
$\|Du\|\leq\|u\|'_{\rho'}<1$.
Hence, $(I+Du)^{-1}=\sum_{n\geq0}(-Du)^n$.
Using the following formulas:
$$
\begin{array}{l}
(Du+Dh)^n=Du^n+\sum\limits_{j=0}^{n-1}Du^j\,Dh\,Du^{n-1-j}+\OO({\|h\|'_{\rho'}}^{2}),
\quad n\geq1,
\\
\sum\limits_{n=1}^{+\infty}\sum\limits_{j=0}^{n-1}(-1)^n\,Du^j\,Dh\,Du^{n-1-j}=-(I+Du)^{-1}\,Dh\,(I+Du)^{-1},
\end{array}
$$
and the Taylor expansion of $f$ around $\theta+u(\theta)$, $\theta\in\DD(\rho)$, that gives
$$
f\circ(\id+u+h)=f\circ(\id+u)+Df\circ(\id+u)\,h+\OO({\|h\|'_{\rho'}}^{2}),
$$
we get
$$
\begin{array}{rl}
F(u+h)-F(u)  
= &
\Ii^-
\sum\limits_{n=0}^{+\infty}(-1)^{n}
\left[(Du+Dh)^n(\omega+f)\circ(\id+u+h)    \right.
\\   
&   \left.
- Du^n(\omega+f)\circ(\id+u)  \right]
\\
= &
\Ii^- (I+Du)^{-1}
\left[(\omega+f)\circ(\id+u+h)  \right.
\\
&   \left.
-(\omega+f)\circ(\id+u)    \right.
\\ 
&   \left.
- Dh\,(I+Du)^{-1}(\omega+f)\circ(\id+u+h)
\right]
\\
&
+\OO({\|h\|'_{\rho'}}^2)
\\
= &
\Ii^- (I+Du)^{-1}
\left[Df\circ(\id+u)\,h  
\right.
\\
&   \left.
-Dh\,(I+Du)^{-1}(\omega+f)\circ(\id+u)
\right] 
\\ 
&
+
\OO({\|h\|'_{\rho'}}^2).
\end{array}
$$
That completes the proof.
\cqd

We want to find a solution for the equation 
\begin{equation}\label{F(u)=0}
F(u)=0.
\end{equation}
For that, 
consider a continuous one-parameter family of maps:
$U_\lambda=\id+u_\lambda$, $\lambda\in [0,1]$, with ``initial'' 
condition $U_0=\id$, i.e. $u_0=0$, such that
$$
F(u_\lambda)=(1-\lambda)F(u_0).
$$
Differentiating the above equation in respect to $\lambda$, we get
\begin{equation}\label{homotopy trick eq}
DF(u_\lambda)\frac{du_\lambda}{d\lambda}=-F(0).
\end{equation}

\begin{remark} \label{remark - invert DF(u) = commutator eq} 
The derivative of $F$ at $u$ can be rewritten in the form
$$
DF(u)\,h=\Ii^-(DU)^{-1}[DX\circ U\,h-Dh\,(DU)^{-1}\,X\circ U].
$$
Evaluating on $v\circ U$, where $v$ is the vector field that
generates $U$, i.e. $dU/d\lambda=v\circ U$ and $v\circ U=DU\,v$, one gets
$$
\begin{array}{ll}
DF(u)\,v\circ U & =\Ii^-(DU)^{-1}\,[v,X]\circ U \\
&=\Ii^-[v,(DU)^{-1}\,X\circ U],
\end{array}
$$
where $[v,w]=Dw\,v-Dv\,w$ is the commutator for vector fields $v,w$.
Thus, writing $\tilde X=(DU)^{-1}\,X\circ U$, it suffices to solve the
commutator equation 
$$
\Ii^-[v,\tilde X]=-\Ii^-X
$$ 
with respect to $v$ satisfying $\Ii^+v=0$, or, equivalently, invert
$DF(u)$.
Note also that we allow to have $v=v_\lambda$.
\end{remark}

\begin{proposition}\label{DF(u) invert}
If $u\in\BB$, then $DF(u)^{-1}$
is a bounded linear operator from
$\Ii^-\A_d(\rho')$ to $\Ii^-\A'_d(\rho')$ and
$$
\|DF(u)^{-1}\| < \frac {\hat\delta}{\hat\varepsilon}.
$$
\end{proposition}

From the above proposition (to be proved in Section \ref{proof of
prop}) we know that $DF(u)$ is invertible for
$u\in\BB$, thus we may
integrate (\ref{homotopy trick eq}) with respect to $\lambda$,
obtaining the integral equation:
\begin{equation}\label{u_l}
u_\lambda=-\int_0^\lambda DF(u_\mu)^{-1}\,F(0)\,d\mu.
\end{equation}
In order to check that $u_\lambda\in\BB$ for any $\lambda\in[0,1]$, we
estimate its norm:
$$
\begin{array}{rl}
\|u_\lambda\|'_{\rho'} 
\leq 
&
\sup\limits_{v\in\BB}\|DF(v)^{-1}F(0)\|'_{\rho'}
\\
\leq&
\sup\limits_{v\in\BB}\|DF(v)^{-1}\|.\|\Ii^-f\|_{\rho'}
\\
< &
\frac{\hat\delta\hat\varepsilon}{\hat\varepsilon}=\hat\delta.
\end{array}
$$
Therefore, the solution of (\ref{F(u)=0}) exists in $\BB$ and is given
by (\ref{u_l}) when $\lambda=1$.

Now, the open ball $\hat B$ as claimed
is simply given by 
$\hat B=\omega +\EE\subset\A'(\rho)$ using $\hat\varepsilon$ given by (\ref{bound on epsilon}).
For $X=\omega+f\in \hat B$ we have $\|u\|'_{\rho'}=\OO(\|f\|_{\rho'})$.
Thus,
\begin{eqnarray*}
\|\UU(X)-\omega\|_{\rho'}
&\leq&
\|\sum_{n\geq1}(-Du)^n\omega+(DU)^{-1}f\circ U\|_{\rho'}\\
&\leq&
\frac1{1-\|u\|'_{\rho'}}
(\|\omega\|\,\|u\|'_{\rho'}+\|f\|_\rho)\\
&\leq&
2\left(
\frac{\hat\delta}{\hat\varepsilon}\|\omega\|\|f\|_{\rho'}+
\|f\|_{\rho'}
\right)\\
&\leq&
2\left(
1+\max\left\{
\frac{2(3-\sqrt6)}{3}, \frac{6\|\omega\|}\sigma (\sqrt6+2)
\right\}
\right)
\|X-\omega\|'_\rho.
\end{eqnarray*}
Notice that $\UU(X)=(DU)^{-1}(\omega+f\circ U)$.
Moreover,
$\|\UU(X)-\omega-\Ii^+f\|_{\rho'}=\OO(\|f\|_{\rho'}^2)$.
This means that the derivative at $\omega$ of the map $\UU$
is $\Ii^+$.
Now, assume $\|f-\psi\|'_\rho=\OO(\|f\|_{\rho'})$, $\omega+\psi\in \hat B\cap\Cc^d$, and write
$\UU(X)=\UU(\omega+\psi+(f-\psi))$.
So, 
$$
\|\UU(\omega+\psi+(f-\psi))-(\omega+\psi)-\Ii^+(f-\psi)\|_{\rho'} =
\OO(\|f\|_{\rho'}^2),
$$
i.e. $\Ii^+$ is the derivative of $\UU$ at $\omega +\psi$.
That completes the proof of Theorem \ref{main theorem1}.

%%%%%%%%%%%%%%%%%%%%%%%%%%%%%%%%%%%%%%%%%%%%%%%%%%%%%%%%%%%%%%%%%%%%%%%%%%%

\subsection{Proof of Proposition \ref{DF(u) invert}} \label{proof of prop}

To prove Proposition \ref{DF(u) invert} we start by inverting the
derivative of $F$ at $u=0$.

\begin{lemma}\label{DF(0) invert}
If $f\in\EE$ such that $\|f\|'_\rho<\sigma/4$,
the associated bounded linear operator $DF(0)^{-1}$, from
$\Ii^-\A_d(\rho')$ to $\Ii^-\A'_d(\rho')$, satisfies:
$$
\|DF(0)^{-1}\| < \frac2{\sigma-4\|f\|'_\rho}.
$$
\end{lemma}

\proof
From Lemma \ref{DF formula} one has 
$$
\begin{array}{rl}
DF(0)\,h
&
=\Ii^-(\hat f\,h-Dh\,\omega)
\\
&=\Ii^- \left[ \hat f-(D\cdot\,\omega)\right]\,h
\\
&=-\left[\Ii-\Ii^-\hat
f\,(D\cdot\,\omega)^{-1}\right](D\cdot\,\omega)\,h,
\end{array}
$$
where $\hat f\,h=Df\,h-Dh\,f$. 
Thus, the inverse of this operator, if it exists, is given by
$$
DF(0)^{-1}=-(D\cdot\,\omega)^{-1}\left[\Ii-\Ii^-\hat
f\,(D\cdot\,\omega)^{-1}\right]^{-1}.
$$
The inverse of $(D\cdot\,\omega)$ is the linear map from
$\Ii^-\A_d(\rho')$ to $\Ii^-\A'_d(\rho')$:
$$
(D\cdot\,\omega)^{-1}\,g(\theta)=\sum\limits_{k\in I^-}\frac{g_k}{2\pi
ik\cdot\omega} e^{2\pi ik\cdot\theta}.
$$
So,
$$
\begin{array}{rl}
\|(D\cdot\,\omega)^{-1}\,g\|'_{\rho'}=
&
\sum\limits_{k\in I^-}\frac{1+2\pi\|k\|}{2\pi
|k\cdot\omega|} \|g_k\|e^{\rho'\|k\|}
\\
<
&
\sum\limits_{k\in I^-}\frac{1+2\pi\|k\|}{2\pi\sigma
\|k\|} \|g_k\|e^{\rho'\|k\|} 
\\
\leq
&
\frac2\sigma\|g\|_{\rho'},
\end{array}
$$
where the use of the definition of $I^-$ (Definition \ref{far from
resonance terms}) was crucial to avoid dealing with arbitrarily small denominators.
Hence, 
$\|(D\cdot\,\omega)^{-1}\| < \frac2\sigma$.

Similarly, for $\hat f\colon
\Ii^-\A'_d(\rho')\to\A_d(\rho')$,
$$
\|\hat f\,h\|_{\rho'}\leq
\|f\|'_{\rho'}\|h\|_{\rho'}+\|f\|_{\rho'}\|h\|'_{\rho'}
$$
Hence it is a bounded operator with $\|\hat f\|\leq 2\|f\|'_{\rho'}$.
Therefore,
$$
\|\Ii^-\hat f(D\cdot\,\omega)^{-1}\| < \frac4\sigma\|f\|'_{\rho'}<1,
$$
assuming $\|f\|'_\rho <\sigma/4$.
So, $\left[\Ii-\Ii^-\hat f(D\cdot\,\omega)^{-1}\right]^{-1}$ is a
bounded linear map in $\Ii^-\A_d(\rho')$ such that:
$$
\left\|\left[\Ii-\Ii^-\hat f(D\cdot\,\omega)^{-1}\right]^{-1}\right\| < \frac\sigma{\sigma-4\|f\|'_{\rho'}},
$$
completing the proof of the Lemma.
\cqd

It remains now to estimate the variation of $DF$ with $u\in\BB$.

\begin{lemma}\label{DF(u)-DF(0) bound}
Let $0<\hat\delta< \min\{\frac12,\frac{\rho-\rho'}{4\pi}\}$.
Given $u\in\BB$ such that $\|u\|'_{\rho'}=\delta<\hat\delta$, 
the linear operator 
$DF(u)-DF(0)$,
mapping $\Ii^-\A'_d(\rho')$ into $\Ii^-\A_d(\rho')$ for $f\in\EE$ with
$\|f\|'_\rho=\varepsilon < \hat\varepsilon$, is bounded and
$$
\|DF(u)-DF(0)\| < \frac{\delta}{1-\delta}
\left[
\left(\frac{4\pi}{\rho-\rho'}+ \frac{4-2\delta}{1-\delta}\right)  \varepsilon
+\frac{2-\delta}{1-\delta}
\|\omega\|
\right].
$$
\end{lemma}

\proof
Lemma \ref{DF formula} gives us that
$$
\begin{array}{rl}
\left[DF(u)-DF(0)\right]\,h= 
&
\Ii^-(I+Du)^{-1}
\left[Df\circ(\id+u)\,h-(I+Du)Df\,h   \right. 
\\
& \left. 
-Dh\,(I+Du)^{-1}(\omega+f)\circ(\id+u) 
\right.
\\
& \left.
+ (I+Du)Dh\,(\omega+f)\right]
\\
= & 
\Ii^-(I+Du)^{-1}\{
\left[Df\circ(\id+u)-Df-Du\,Df\right]\,h
\\
&
+ Du\,Dh\,(\omega+f)
\\
& - Dh\,(I+Du)^{-1}\left[f\circ(\id+u)-f-Du\,(\omega+f)\right]
\}
\\
=&
\Ii^-(I+Du)^{-1}\{A+B+C\},
\end{array}
$$
where $A$, $B$ and $C$ are each of the respective three terms in the
previous sum.
Thus,
$\|(I+Du)^{-1}\|\leq(1-\|Du\|)^{-1}$.
Using Lemma \ref{lemma f-f*}, one obtains:
$$
\|A\|_{\rho'} \leq
\left( \frac{4\pi}{\rho-\rho'}\|f\|'_{\rho}\|u\|_{\rho'}+  \|f\|'_{\rho'}\|u\|'_{\rho'} \right)\|h\|_{\rho'}
$$
and $\|B\|_{\rho'} \leq \left(\|\omega\|+\|f\|_{\rho'}\right) \|u\|'_{\rho'}\|h\|'_{\rho'}$.
Finally, noting that
$$
\|Du\,(I+Du)^{-1}\|\leq \|Du\|.(1-\|Du\|)^{-1},
$$
and from Lemma \ref{lemma f-f*} again, one gets
the bound for the third term: 
$$
\|C\|_{\rho'} \leq
\frac{1}{1-\|u\|'_{\rho'}}   
\left[
\|f\|'_{(\rho+\rho')/2}\|u\|_{\rho'}+\|u\|'_{\rho'}
\left( \|\omega\|+\|f\|_{\rho'} \right) \right]
  \|h\|'_{\rho'}.
$$
\cqd

To conclude the proof of Proposition
\ref{DF(u) invert}, notice that,
for $u\in\BB$ and $f\in\EE$:
$$
\begin{array}{ll}
\|DF(u)^{-1}\| 
&
\leq \left( \|DF(0)^{-1}\|^{-1}-\|DF(u)-DF(0)\|\right)^{-1}
\\
&
< \left\{ 

\frac\sigma2-2\hat\varepsilon-
\frac{\hat\delta}{1-\hat\delta}
\left[
\left(\frac{4\pi}{\rho-\rho'}+ \frac{4-2\hat\delta}{1-\hat\delta}\right)  \hat\varepsilon
+\frac{2-\hat\delta}{1-\hat\delta}
\|\omega\|
\right]
\right\}^{-1}
\\
&
< \frac{\hat\delta}{\hat\varepsilon}.
\end{array}
$$
The last inequality is true if we choose $\hat\varepsilon$ and $\hat\delta$ to satisfy
\begin{equation}\label{2nd bound on epsilon}
\hat\varepsilon<  
\hat\delta\left[
\frac\sigma{2}-\frac{2\hat\delta}{(1-\hat\delta)^2}\|\omega\| \right]
\left[ 1+2\hat\delta+\frac{\hat\delta^2}{1-\hat\delta} \left( \frac{4\pi}{\rho-\rho'}+\frac{4-2\hat\delta}{1-\hat\delta}\right)
\right]^{-1},
\end{equation}
and
$$
\hat\delta< 
\frac12- \frac12\left[ 1 + \frac\sigma{2\|\omega\|}\right]^{-\frac12} < \frac12.
$$

An effective value of $\hat\varepsilon$ for each $\sigma$ can then be determined by the minimum
value of the upper bounds given in Lemma \ref{DF(0) invert} and (\ref{2nd bound on epsilon}).
That can be obtained from a specific choice of
$\hat\delta$ given by 
$$
\hat\delta = \min \left\{ 
\frac{\rho-\rho'}{4\pi} , \frac{3-\sqrt6}6 \frac {\sigma}{\|\omega\|}
\right\} < \frac12- \frac12\left[ 1 +
\frac\sigma{2\|\omega\|}\right]^{-\frac12} < \frac12,
$$
as long as $\sigma/\|\omega\|<1$.
Therefore, it is sufficient to impose
\begin{equation}\label{bound on epsilon}
\hat\varepsilon:=
\min\left\{
\frac\sigma4,
\frac{\sqrt6-2}{12}\sigma\hat\delta
\right\},
\end{equation}
where the following inequalities were applied to (\ref{2nd bound on epsilon}):
$$
\frac1{1-\hat\delta} < 
\sqrt{1+\frac{\sigma}{2\|\omega\|}}
\quad\text{ and }\quad
\frac{2\hat\delta}{(1-\hat\delta)^2}\|\omega\|<\frac{(3-\sqrt6)\sigma}{2}
<\frac\sigma2.
$$

%%%%%%%%%%%%%%%%%%%%%%%%%%%%%%%%%%%%%%%%%%%%%%%%%%%%%%%%%%%%%%%%%%%%%%%%%%%%%

\subsection{Analyticity Improvement} \label{Analyticity Improvement}

We want to prove Proposition \ref{proposition T}, i.e. that $X\circ T$ is analytic in
$\DD(\rho)$ and has bounded derivative.
Let $\sigma>0$ such that
$$
\max\left\{\frac{\|T^*x\|}{\|x\|} \colon\|x\|\not=0, |\omega\cdot
x|\leq\sigma\|x\| \right\} < \kappa.
$$
Using $X\circ T(\theta)=\sum\limits_{k\in I_\sigma^+}X_k e^{2\pi iT^*k\cdot
\theta}$, 
$$
\begin{array}{rl}
\|X\circ T\|_\rho  & \leq  \sum\limits_{k\in I_\sigma^+}\| X_k\|e^{\rho\|T^*k\|} 
\\
&
\leq \sum\limits_{k\in I_\sigma^+}\|
X_k\|e^{(\rho\kappa-\rho')\|k\|}e^{\rho'\|k\|} 
\leq \|X\|_{\rho'},
\end{array}
$$
and for the derivative we have
$$
\begin{array}{rl}
\|D(X\circ T)\|_\rho  & \leq 2\pi \sum\limits_{k\in I_\sigma^+}\|T^*k\|
e^{-\beta\|k\|} \|
X_k\|e^{(\rho\kappa+\beta-\rho')\|k\|}e^{\rho'\|k\|} 
\\
&
\leq \frac{2\pi}{\beta}\kappa\|X\|_{\rho'},
\end{array}
$$
by choosing $0<\beta<\rho'-\kappa\rho$, and using the relation $\sup_{t\geq0}te^{-\xi t}\leq 1/\xi$ for $\xi>0$.
These bounds imply that
$\| \TT(X) \|'_\rho \leq
(1+{2\pi\kappa}/\beta)\|T^{-1}\|\,\|X\|_{\rho'}$.

Let $r>\rho$ such that $\rho'>r\kappa$. 
What was done above for $\DD(\rho)$ applies as well to $\DD(r)$.
Therefore, one can decompose $\TT=\II\circ\JJ$, where
$\JJ\colon\Ii^+\A(\rho')\to\A'(r)$ as before, and
$\II\colon\A'(r)\to\A'(\rho)$ is the inclusion map
$\II(X)=X|_{\DD(\rho)}$. Note that $\JJ$ is bounded and $\II$ is
compact,
thus completing the proof.

%%%%%%%%%%%%%%%%%%%%%%%%%%%%%%%%%%%%%%%%%%%%%%%%%%%%%%%%%%%%%%%%%%%%%%%%%%%%%

\section{Spectral Properties of $D\RR(\omega)$ and Invariant Manifolds}\label{Spectral Properties and Invariant Manifolds}

As $\sigma$ remains fixed, we continue to drop the index of the sets
$I^\pm$ and the projections $\Ii^\pm$.

First we study the behaviour of the linear map $\LL$ whose spectral
properties are closely related to the ones of $D\RR(\omega)$.
Notice that $\Ee\,\A'(\rho)=\Cc^d$ and $(\Ii-\Ee\,)\A'(\rho)$ are invariant subspaces of $\LL$ as $\LL\circ\Ee=\Ee\circ\LL$.

\begin{lemma}\label{lemma bound on stable part}
If $\omega\in DC(\beta)$, for any $\beta>0$, and $\sigma\|\bar\omega\|<\frac12$,
one finds constants $a>0$ and $b,c>1$ such that
$$
\|\LL^{n+1}(\Ii-\Ee)\| \leq b \lambda_1^2 e^{- a c^n}
\|\LL^{n}(\Ii-\Ee)\|, \quad n\geq0
$$
\end{lemma}

\proof
For every $n\geq0$, we define
$I^+_n$ as the set of indices $k$ that verify $(T^*)^{n}k\in
I^+\setminus\{0\}$.
Hence,
$$
\sigma\|(T^*)^{n}k\| \geq
|\omega\cdot (T^*)^{n}k| =
|T^{n}\omega\cdot k|=
|\lambda_1|^{n}|\omega\cdot k|, \quad \text{ for } k \in I^+_n.
$$
This inequality implies that
$$
|\omega\cdot k| \leq \sigma \|(\omega\cdot k)\,\bar\omega\| +
\sigma\norm { \left(\lambda_1^{-1}T^*\right)^n[k-(\omega\cdot k)\,\bar\omega]}.
$$
Thus, assuming $\sigma\|\bar\omega\|<\frac12$ and knowing that
$k-(\omega\cdot k)\,\bar\omega$ is the component of $k\in I^+_n$ on the
spectral directions of $T^*$ corresponding to the eigenvalues
$\lambda_j$, $j=2,\dots,d$,
$$
|\omega\cdot k| \leq
\frac{\sigma}{1-\sigma\|\bar\omega\|}\abs{\frac{\lambda_2}{\lambda_1}}^n
\|k-(\omega\cdot k)\,\bar\omega\|
$$
As $\omega$ is a diophantine vector of degree $\beta>0$, there is
a constant $C>0$ such that
$$
|\omega\cdot k| > C \|k\|^{-\beta-1},
$$ 
and we get a lower bound for the norm of $k\in I^+_n$,
\begin{equation}\label{lower bound for k in I+n}
\|k\|^{2+\beta} > C \frac{1-\sigma\|\bar\omega\|}{2\sigma}
\abs{\frac{\lambda_1}{\lambda_2}}^n.
\end{equation}

\begin{figure}
\begin{center}
\input {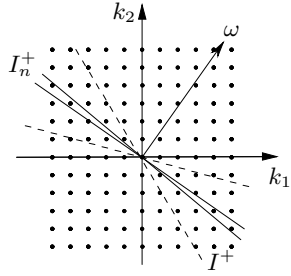} 
\end{center}
\caption{The set $I_n^+\subset I^+$ for the two-dimensional case.}     
\label{I_n+}
\end{figure}

Let us now define the operator $\Ii^+_n\colon \A'(r)\to\A'(\rho)$,
where we have chosen $r>\rho$. 
This is a projection for the indices in $I^+_n$ together
with an analytic inclusion. 
Making use of the bound (\ref{lower bound for k in I+n}),
the operator norm of $\Ii^+_n$ follows from 
$$
\|\Ii^+_{n}f\|'_\rho = \sum\limits_{I^+_{n}}(1+2\pi\|k\|)\|f_k\|
e^{r\|k\|}e^{-(r-\rho)\|k\|} 
\leq e^{-ac^n} \|f\|'_r,
$$
with
$a=(r-\rho)[C(2\sigma)^{-1}(1-\sigma\|\bar\omega\|)]^{1/(2+\beta)}>0$
and $c=|\lambda_1\lambda_2^{-1}|^{1/(2+\beta)}>1$.

We can write $\LL^{n+1}$ with respect to $\LL^n$, including
the operator $\Ii^+_n$, in the form
$$
\LL^{n+1}(\Ii-\Ee)f=\LL^n(\Ii-\Ee)\Ii^+_{n}\bar\LL f,
$$
where $\bar\LL\colon\A'(\rho)\to\A'(r)$ is $\LL$ followed
by an analytic extension.
That is, $\bar\LL=\lambda_1\bar\TT\circ\Ii^+$ where
$\bar\TT\colon\Ii^+\A'(\rho)\subset\Ii^+\A(\rho)\to\A'(r)$ given as in
Proposition \ref{proposition T}.
The norm $\|\bar\LL\|\leq \|\bar\TT\|$ is determined by
$$
\|\bar\TT f\|'_r\leq b \|T^{-1}\| \|f\|'_\rho,
$$
for some $r>\rho$ and a constant $b>1$ from the proof of Proposition
\ref{proposition T} in Section \ref{Analyticity Improvement}.
So, the claim follows from
$$
\|\LL^{n+1}(\Ii-\Ee)f\|'_\rho  \leq
\|\LL^n(\Ii-\Ee)\|\,\|\Ii^+_{n}\|\, \|\bar\LL\| \,\|f\|'_\rho,
$$
and $\|T^{-1}\| \leq |\lambda_1|$.
\cqd

So, the direction given by the non-constant terms $\Ii-\Ee$ is a
stable eigenspace with eigenvalue zero.
The spectrum of $\LL\circ \Ee$ is simply the one of $T^{-1}$.
That is, the eigenvalues are $\lambda_1\lambda_j^{-1}$, $j=2,\dots,d$, and
$1$ corresponding to the eigenvector $\omega$.
Hence, there is a ``neutral'' direction given by $\omega$ whereas the
remaining ones corresponding to $\omega^{(2)},\dots,\omega^{(d)}$ are all unstable.

The spectral properties of $D\RR(\omega)$ are easily related to those
of $\LL$ because
$$
D\RR(\omega)=(\Ii-\Pp_\omega\Ee)\LL,
$$
where $\Pp_\omega$ is the projection on the subspace spanned by
$\omega$, i.e. $\Pp_\omega f_0=(\omega\cdot f_0)\,\bar\omega$.
Note also that $\LL$ commutes with the projection $\Pp_\omega\Ee$
because $\omega$ is a eigenvector of $T$, thus of $\LL$. So,
$$
(\Ii-\Pp_\omega\Ee)\LL = \LL (\Ii-\Pp_\omega\Ee).
$$

Therefore, the projection $\Pp^u=(\Ii-\Pp_\omega)\Ee$ determines 
the linear space $\Pp^u\A'(\rho)$ spanned by
$\{\omega^{(2)},\dots\omega^{(j)}\}$, which is the unstable eigenspace
of the linearised map at $\omega$.
In a similar way, $\Pp^s=\Ii-\Pp^u$ is the stable part.

We are now in a position to define and determine the stable and unstable invariant
local manifolds of $\RR$ at the fixed point $\omega$. 
The stable one, $\WW^s(\omega)$, corresponds to all the vector fields
$X$ in some neighbourhood $B$ such that
$\RR^n(X)\in B$, for any $n\geq0$.
Similarly, considering the inverse iterates of $\RR$, we define the 
unstable local manifold $\WW^u(\omega)$.

If $X\in\WW^s(\omega)$, the winding set is equal to that of $\omega$,
since it is preserved by $\RR$. 
This is in fact true because $\omega$
is an eigenvector of the matrix $T\in GL(d,\Zz)$, the linear part of the
transformation induced by $\RR$ on the domain of $X$. 
The winding
set is preserved up to the action of $\hat T\colon x\mapsto Tx/\|Tx\|$, by iterating $\RR$.
If the winding set of a vector $X\in B$ is $w_X=\omega/\|\omega\|$,
then it is also in $\WW^s(\omega)$, otherwise would approximate the
unstable direction and $w_X$ would be different.

The set of constant vector fields is invariant under $\RR$,
i.e. $\RR(\Cc^d \cap B) \subset \Cc^d$, as can be seen
from the
way the non-linear coordinate change $\UU$ was constructed in Theorem
\ref{main theorem1}.
This allows us to conclude that the invariant unstable local manifold is
$\WW^u(\omega) = \Pp^u (B)$.

\begin{lemma}
We can find a ball $B'\subset \A'_d(\rho)$ around $\omega$ such that,
for any $X\in\WW^s(\omega)\cap B'$, there are constants $K>0$ and $\theta\in]0,1[$ yielding
$$
\|\RR^m(X)-\omega\|'_\rho < K \theta^m \|X-\omega\|'_\rho, \quad m>0.
$$
Moreover, there exists $a>0$, $c>1$ and $N>0$ satisfying
$$
\| (\Ii-\Ee) \RR^n(X)\|'_\rho < e^{-a c^n} \|X-\omega\|'_\rho + \theta^{2n} {\|X-\omega\|'_\rho}^2,\quad n>N.
$$
\end{lemma}

\proof
Fixing $0<t'<t<\theta<1$ and choosing an appropriate norm $\|\cdot\|$ on $\A'(\rho)$ equivalent to $\|\cdot\|'_\rho$, it is possible to have $\|D\RR(\omega)\, f\| < t' \|f\|$, $f\in\Pp^s(\A'(\rho))$, since $D\RR(\omega)$ is compact.

We can write $X=\omega+f+g$ where $f=\Pp^s(X-\omega)$ and $g=\Pp^u(X-\omega)$.
By the stable manifold theorem, there is a constant $A>0$ for which 
\begin{equation}\label{bounds on g and f}
\|g\| < A\,{\|f+g\|}^2,\quad \text{ and } \quad
\|f\|< \|f+g\| + A \|f+g\|^2.
\end{equation} 
The image of $X$ is then $X'=\RR(X)=\omega+f'+g'$ with
\begin{equation}\label{f' and g'}
\begin{array}{l}
f'=D\RR(\omega)\, f,  \\
g'=D\RR(\omega)\, g + \OO({\|f+g\|}^2).
\end{array}
\end{equation}
The analyticity of $\RR$ and the formula of the second-order Taylor
remaining imply that there exists constants $B,C>0$ satisfying
\begin{equation}\label{bounds on f' and g'}
\begin{array}{l}
\|f'\| < t' \|f\|,  \\
\|g'\| < C \, \|g\| + B \, {\|f+g\|}^2.
\end{array}
\end{equation}
By combining this with (\ref{bounds on g and f}) we obtain 
$$
\|f'+g'\|<t'\|f+g\| + (At'+AC+B) \|f+g\|^2.
$$
So, there is a radius $r>0$ such that 
$$
\|X'-\omega\|< t \|X-\omega\|
$$ 
for $\|X-\omega\|<r$.
Iterating this inequality $m$ times we obtain
$$
\|\RR^m(X)-\omega\|'_\rho < K \theta^m \|X-\omega\|'_\rho,
$$
if $B'=\{X\in \A'_d(\rho) \colon \|X-\omega\|<r\}$, that is the first of our claims.

From Lemma \ref{lemma bound on stable part} we can find $a>0$ and $c>1$
and a sequence of positive integers $k_n<\varepsilon n$, $\varepsilon>0$, such that
$$
\|(\Ii-\Ee)[D\RR(\omega)]^{k_n}\| < e^{-ac^n},
$$
whenever $n$ is sufficiently large.
Given such $n$ and $X\in \WW^s(\omega) \cap B'$,
the first part of the proof showed that $\|\RR^m(X)-\omega\| <\delta = t^m r$, with $m=n-k_n$.
Consider a map $F\colon (f, g)\mapsto (f',g')$ defined by
(\ref{f' and g'}), and let $f_0=\Pp^s(\RR^m(X)-\omega)$ and $g_0=\Pp^u(\RR^m(X)-\omega)$.
The pair $(f_0,g_0)$ satisfies 
$$
\begin{array}{l}
\|f_0\|<2\delta, \\
\|g_0\|< A \delta^2.
\end{array}
$$
Applying the inequalities
(\ref{bounds on f' and g'}),
if $1<A<1/\delta$, then $\|g_1\| < (C+9B)A\delta^2$, using the notation $(f_i,g_i)=F^i(f_0,g_0)$, $i\geq0$.
The fact that $\|f_i\|$ is bounded by at least $2\delta$ allows us to 
iterate $F$ for $k_n$ times, obtaining: 
$$
\begin{array}{l}
f_{k_n} =  [D\RR(\omega)]^{k_n} f_0,  \\
\|g_{k_n}\| < (C+9B)^{k_n} A \delta^2,
\end{array}
$$
as long as $(C+9B)^{k_n}A\delta<1$. This is in fact verified if we choose $\varepsilon>0$ such that $(C+9B)^{k_n} t^{2m} < \theta^{2n}$, or, more strongly, 
$$
[t^{-2}(C+9B)]^\varepsilon < (\theta/ t)^2.
$$

Now, $\RR^n(X) = \omega + f_{k_n} + g_{k_n}$ and $(\Ii-\Ee)\RR^n(X) = (\Ii-\Ee)(f_{k_n} + g_{k_n})$. Hence,
$$
\|(\Ii-\Ee)\RR^n(X)\| < 2 e^{-ac^n} t^m r + A\theta^{2n} r^2,
$$
which completes the proof, if $n$ is chosen sufficiently large.
\cqd

%%%%%%%%%%%%%%%%%%%%%%%%%%%%%%%%%%%%%%%%%%%%%%%%%%%%%%%%%%%%%%%%%%%%%%%%%%%%
\section*{Acknowledgements} 

I would like to express my gratitude to Professor R. S. MacKay for 
the orientation and the support given. 
I also wish to thank Cristel Chandre and Tim Hunt for
useful conversations.
The author is funded by Funda\c c\~ao para a Ci\^encia e
a Tecnologia, under the research grant BD/11230/97.
%%%%%%%%%%%%%%%%%%%%%%%%%%%%%%%%%%%%%%%%%%%%%%%%%%%%%%%%%%%%%%%%%%%%%%%%%%%%

\bibliographystyle{plain} 
\addcontentsline{toc}{section}{References}
\bibliography{rfrncs}

\begin{thebibliography}{1}

\bibitem{Arnold2}
V.~I. Arnol'd.
\newblock Small denominators {I}, mappings of the circumference onto itself.
\newblock {\em Transl. {AMS} 2nd Series}, 46:213--284, 1961.

\bibitem{BGKM}
C.~Baesens, J.~Guckenheimer, S.~Kim, and R.~S. MacKay.
\newblock Three coupled oscillators: mode-locking, global bifurcations and
  toroidal chaos.
\newblock {\em Physica D}, 49:387--475, 1991.

\bibitem{Delshams}
A.~Delshams and R.~de~la Llave.
\newblock {KAM} theory and a partial justification of {G}reene's criterion for
  non-twist maps.
\newblock {\em SIAM J. Math. Anal.}, 31(6):1235--1269, 2000.

\bibitem{Katok}
A.~Katok and B.~Hasselblatt.
\newblock {\em Introduction to the modern theory of Dynamical Systems}.
\newblock Cambridge University Press, 1995.

\bibitem{Koch}
H.~Koch.
\newblock A renormalization group for {H}amiltonians, with applications to
  {KAM} tori.
\newblock {\em Erg. Theor. Dyn. Syst.}, 19:475--521, 1999.

\bibitem{Lang}
S.~Lang.
\newblock {\em Introduction to diophantine approximations}.
\newblock Springer-Verlag, 2nd edition, 1995.

\bibitem{jld2}
J.~Lopes~Dias.
\newblock Renormalisation of vector fields for a generic frequency vector.
\newblock {\em preprint http://arXiv.org/abs/math.DS/0105067}, 2001.

\bibitem{MacKay}
R.~S. MacKay.
\newblock Three topics in {H}amiltonian dynamics.
\newblock In Y.~Aizawa, S.~Saito, and K.~Shiraiwa, editors, {\em Dynamical
  Systems and Chaos}, volume~2. World Scientific, 1995.

\bibitem{Poschel}
J.~P\"oschel.
\newblock On elliptic lower dimensional tori in {H}amiltonian systems.
\newblock {\em Math. Z.}, 202:559--608, 1989.

\end{thebibliography}

\end{document}